\documentclass[11pt]{article}
\usepackage[english]{babel} 
\usepackage[utf8]{inputenc} 
\usepackage[a4paper,left=2.5cm,bottom=2.5cm,top=2cm]{geometry} 
\usepackage{amsmath} 
\usepackage{amssymb}
\usepackage{amsfonts}
\usepackage{bm} 
\usepackage{graphicx} 
\usepackage[export]{adjustbox} 
\usepackage{pdflscape} 
\usepackage{fancyhdr} 
\usepackage{color} 
\usepackage{wrapfig} 
\usepackage{lipsum} 
\usepackage[normalem]{ulem} 
\usepackage{blindtext}
\usepackage{titlesec}
\usepackage{algpseudocode}
\usepackage{enumitem}
\usepackage{caption}
\usepackage{subcaption}
\usepackage{tcolorbox, todonotes}
\usepackage{glossaries}
\usepackage{csvsimple}
\usepackage{multicol}
\usepackage{multirow}
\usepackage{adjustbox}
\usepackage{booktabs}
\usepackage[bookmarks=false]{hyperref}
\usepackage[toc,page]{appendix}

\usepackage{colortbl}
\usepackage{appendix}
\usepackage{float}
\usepackage[numbers, square, sort&compress]{natbib}
\usepackage{url}

\usepackage{authblk}

\DeclareMathOperator*{\argmax}{arg\,max}

\usepackage[symbol]{footmisc}
\newtheorem{theorem}{Theorem}
\newtheorem{proposition}[theorem]{Proposition}

\newtheorem{definition}{Definition}

\newenvironment{pf}{\noindent{\bf Proof. }}{\hfill $\blacksquare$ \\}

\RequirePackage{etoolbox}
\AtBeginEnvironment{tabular}{
\sffamily\fontsize{7.5}{10}\selectfont
}

\allowdisplaybreaks

\newcommand{\states}{\mathcal{S}}
\newcommand{\state}{h}
\newcommand{\action}{a}
\newcommand{\actions}{\mathcal{A}}

\newcommand{\timeset}{\mathcal{T}}
\newcommand{\epoch}{t}
\newcommand{\transprobs}{P}

\newcommand{\rewards}{\mathcal{R}}

\newcommand{\discount}{\alpha}
\newcommand{\grayarea}{\Lambda}
\newcommand{\obs}{\mathcal{N}}
\newcommand{\val}{v}

\usepackage{tikz}
\usetikzlibrary{decorations.pathreplacing,calc,arrows}

\definecolor{ForestGreen}{RGB}{34,139,34}

\newacronym[longplural=Markov decision processes]{mdp}{MDP}{Markov decision process}
\newacronym[longplural = multi-model MDPs] {mmdp}{MMDP}{multi-model Markov decision process}
\newacronym{cmdp}{CMDP}{contextual Markov decision process}
\newacronym[longplural = partially-observable MDPs]{pomdp}{POMDP}{partially-observable MDP}
\newacronym[longplural = decision makers]{dm}{DM}{decision maker}
\newacronym[longplural = transition probability matrices]{tpm}{TPM}{transition probability matrix}
\newacronym[longplural = discrete time Markov chains]{dtmc}{DTMC}{discrete time Markov chain}
\newacronym{rl}{RL}{reinforcement learning}

\newacronym{mip}{MIP}{mixed-integer program}
\newacronym{lp}{LP}{linear program}
\newacronym{ws}{WS}{wait-and-see}
\newacronym{dfr}{DFR}{Decreasing Failure Rate}
\newacronym{adp}{ADP}{approximate dynamic programming}

\title{The Implications of State Aggregation in Deteriorating Markov Decision Processes with Optimal Threshold Policies}

\author[a]{Madeleine Pollack}
\author[b]{Lauren N. Steimle}

\affil[a]{Massachusetts Institute of Technology\protect \\ {\small \tt \{pollack9\}@mit.edu}}
\affil[b]{Georgia Institute of Technology\protect \\ {\small \tt \{steimle\}@gatech.edu}}

\date{\today}

\begin{document}

\maketitle

\begin{abstract}
    Markov Decision Processes (MDPs) are mathematical models of sequential decision-making under uncertainty that have found applications in healthcare, manufacturing, logistics, and others. In these models, a decision-maker observes the state of a stochastic process and determines which action to take with the goal of maximizing the expected total discounted rewards received. In many applications, the state space of the true system is large and there may be limited observations out of certain states to estimate the transition probability matrix. To overcome this, modelers will aggregate the true states into ``superstates" resulting in a smaller state space. This aggregation process improves computational tractability and increases the number of observations among superstates. Thus, the modeler’s choice of state space leads to a trade-off in transition probability estimates. While coarser discretization of the state space gives more observations in each state to estimate the transition probability matrix, this comes at the cost of precision in the state characterization and resulting policy recommendations. In this paper, we consider the implications of this modeling decision on the resulting policies from MDPs for which the true model is expected to have a threshold policy that is optimal. We analyze these MDPs and provide conditions under which the aggregated MDP will also have an optimal threshold policy. Using a simulation study, we explore the trade-offs between more fine and more coarse aggregation. We explore the the show that there is the highest potential for policy improvement on larger state spaces, but that aggregated MDPs are preferable under limited data. We discuss how these findings the implications of our findings for modelers who must select which state space design to use.
\end{abstract}

\pagenumbering{arabic}

Markov Decision Processes (MDPs) are mathematical models of sequential decision-making under uncertainty that have found applications in healthcare, manufacturing, logistics, and others. In these models, a decision-maker observes the state of a stochastic process and determines which action to take with the goal of maximizing the expected total discounted rewards received. In many applications, the state space of the true system is large and there may be limited observations out of certain states to estimate the transition probability matrix. To overcome this, modelers will aggregate the true states into ``superstates" resulting in a smaller state space. This aggregation process improves computational tractability and increases the number of observations among superstates. Thus, the modeler’s choice of state space leads to a trade-off in transition probability estimates. While coarser discretization of the state space gives more observations in each state to estimate the transition probability matrix, this comes at the cost of precision in the state characterization and resulting policy recommendations. In this paper, we consider the implications of this modeling decision on the resulting policies from MDPs for which the true model is expected to have a threshold policy that is optimal. We analyze these MDPs and provide conditions under which the aggregated MDP will also have an optimal threshold policy. Using a simulation study, we explore the trade-offs between more fine and more coarse aggregation. We explore the the show that there is the highest potential for policy improvement on larger state spaces, but that aggregated MDPs are preferable under limited data. We discuss how these findings the implications of our findings for modelers who must select which state space design to use.

\section{Background} \label{section: problem statement}

In this section, we describe discrete-time infinite-horizon MDPs and explain the method for deriving the estimated TPM for a given MDP model of this type. We also give an introduction on how to aggregate the state space of a Markov chain or MDP model.

\subsection{Markov decision processes}\label{section: mdp}
A Markov decision process (MDP) is a model of a stochastic control process in which the DM seeks to take actions to control a stochastic system. The DM observes the system at discrete time points, $\timeset = \{0,1,2,3,\ldots\}$  where $\epoch \in \timeset$ represents the \textit{decision epoch} or amount of time (e.g., months, years, etc.) that has passed since the beginning of the planning horizon. 
In this article, we focus on infinite-horizon MDP models with a countably infinite number of decision epochs. At each decision epoch, the \textit{health state}, or simply, \textit{state} of the system $\state \in \states$ is observed, where the \textit{state space}, $\states$, is the set of all possible states. After observing the state of the system, the DM takes some action $a \in \actions$, where the \textit{action space}, $\actions$, is the set of all possible actions. 
Once this action $a$ is taken, the DM receives a real-valued \textit{reward} $r(\state,a)$. The $|\states| \times |\actions|$ reward matrix $\rewards$ contains all possible rewards.
Rewards are discounted at a rate of $\discount \in (0, 1)$ to reflect that rewards received in the future are worth less than those received in the present.
Given the current state $\state$ and the action $a$, the conditional probability $p_{\state \state'}^a \in [0,1]$, often denoted in other literature as  $P(\state' \, | \, \state, a)$, describes the likelihood that the system transitions to a new state $\state' \in \states$ in the next decision epoch. The $|\states| \times |\states| \times |\actions|$ matrix describing the stochastic progression of the system 
makes up the TPM, denoted by $P$. 
We will exclusively consider stationary rewards and transition probabilities, which means that the rewards and transition probabilities are dependent only on the state and action and not on the decision epoch.

For a given realization of a sequence of observed states and subsequent actions
$\bigl((\state_0,a_0), \\ (\state_1,a_1), (\state_2, a_2),...\bigr)$, the realized discounted total reward is given by
\begin{equation}\label{eq. realized reward}
    \sum_{t = 0}^\infty \rewards(\state_t, a_t) = \sum_{t=0}^\infty \discount^t r(\state_t, a_t).
\end{equation}
The goal of the DM is to determine the policy that maximizes the expected value of \eqref{eq. realized reward}.  Because the rewards and transition probabilities are stationary, there will exist an optimal policy that is stationary (i.e., independent of time) and deterministic (i.e. the optimal policy maps each state in the vector to a single optimal action with probability 1) \citep[Theorem 6.2.10]{puterman2014markov}. 

To determine the optimal policy, one must solve the following optimality equations: 
\begin{equation}\label{eq. value to go}
    v^*(\state) = \max_{a \in \actions}\left\{r(\state,a) + \discount \sum_{\state' \in \states} p_{\state\state'}^a v^*(\state')\right\}, \forall \state \in \states.
\end{equation}
Here, the optimal value function $v^*: \state \mapsto \mathbb{R}$ maps a state $\state\in\states$ to the maximum value that the DM can receive if the system starts in state $\state$. The optimal policy $\pi^*$ is a vector of length $|\states|$ with the $\state^\text{th}$ entry containing the action $a$ which maximizes \eqref{eq. value to go} given the system is currently in state $\state$. The optimal policy is given by
\begin{equation}\label{eq. optimal policy}
    \pi^{*}(\state) = \argmax_{a \in \actions} \left\{r(\state, a) + \discount \sum_{h' \in \states}p_{\state\state'}^a v^*(\state')\right\}, \forall \state\in\states.
\end{equation}
Policy iteration, value iteration, and linear programming are common methods of solving  \eqref{eq. value to go} and subsequently \eqref{eq. optimal policy} \citep{puterman2014markov}. 

\subsection{Special MDPs of Interest: 
 Optimal Stopping Time Problems}\label{subsection:mdps of interest}
In this work, we focus on \glspl{mdp} with two main characteristics. The first characteristic is \glspl{mdp} with \glspl{tpm} that are \textit{deteriorating}. This is often called the \gls{dfr} property if the state space is ordered from the least desirable states to the most desirable states. In this paper, our state spaces are ordered such that the \gls{dfr} property holds to remain consistent with the experiments performed by \cite{regnier2013state}. Formally, we can define the \gls{dfr} property as follows.
\begin{definition}[\textbf{Decreasing Failure Rate Property} \citep{regnier2013state}] 
$P^\action$ is said to have the DFR \\ property if its rows are in decreasing stochastic order. 
That is, $P^a$ with state space $\states_J = \{0, 1, 2, \ldots, J\}$ is DFR if $\sum_{j=0}^n p_{ij}^a$ is nonincreasing in $i$ for all $n\in \{0,\ldots,J\}$.
\end{definition}
\noindent Deteriorating \glspl{tpm} are common for application areas like chronic disease progressions \citep{alagoz2004optimal} and probabilities of equipment failures \citep{rahim1993generalized}.

The second characteristic of interest is \glspl{mdp} whose optimal policies are \textit{threshold policies}. This often arises in \textit{stopping time problems} with a deteriorating Markov chain and action space $\{a_1, \;a_2\}$ where $a_1$ is the ``do-nothing" action, and $a_2$ is the ``intervention" action that stops the process. Under certain conditions on the TPM and rewards matrix \cite[\S 6.11.2]{puterman2014markov}, there is guaranteed to exist an optimal policy that is a \textit{threshold policy} (also referred to as a \textit{control-limit policy}), defined as follows.

\begin{definition}[\textbf{Threshold policy} ]\label{def. threshold} 
Given an MDP with ordered state space $\states$ and action space $ \actions = \{a_1,a_2\}$, a \textbf{threshold policy} with \textbf{threshold} $h^*$ is a deterministic Markov policy with an optimal policy $\pi^*$ of the form
\begin{equation*}
    \pi^*(\state)=\begin{cases}
        a_1 \quad &\text{if } \state \leq \state^*\\
        a_2 \quad &\text{if } \state > \state^*.\\
        \end{cases}
        \end{equation*}
\end{definition}
Threshold policies are often assumed to be optimal in \glspl{mdp} with deteriorating \glspl{tpm} because of the intuitive natures of these policies (e.g., in medical decision-making, once a patient reaches a critical level of health, he or she is more likely to continue to deteriorate, so we should initiate treatment). Furthermore, threshold policies are highly interpretable and easily implementable, which makes them attractive to \glspl{dm}. 

\section{The Implications of State Aggregation on MDPs}\label{section: state aggregation}

\begin{figure}[t]
\centering
\includegraphics*[width=.75\textwidth,height= .6\textheight,keepaspectratio]{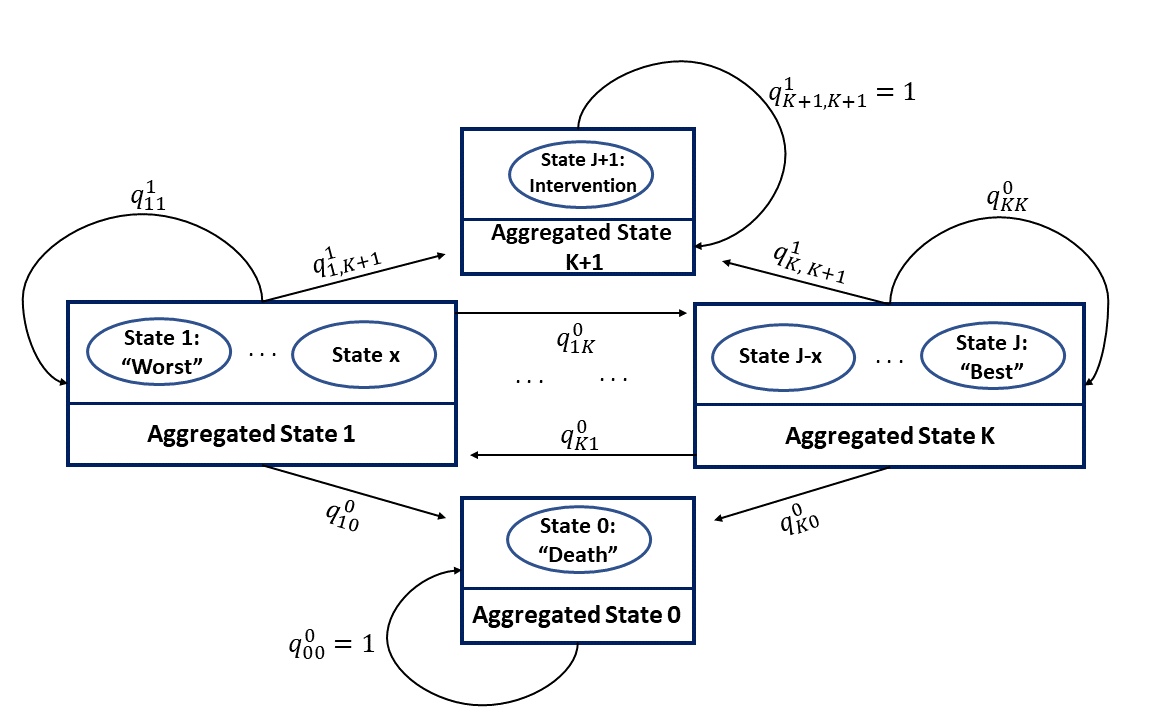}
\vspace{-0.2cm}
\caption{Diagram of aggregated Transition Probability Matrix.\label{fig: lumped TPM}}
\end{figure}

In this section, we formalize the state aggregation process and analyze the properties of aggregated \glspl{mdp}. State aggregation is the process of reducing the state space size by aggregating similar states according to a rule. For ordered state spaces, similar states are typically defined as some number of consecutive states on the ordered state space. An example of an optimal stopping time problem with state aggregation is shown in Figure \ref{fig: lumped TPM}. A valid aggregation of the states into superstates is defined below:

\begin{definition}[\textbf{Valid state aggregation function}]\label{def: agg function} Consider a state space $\states_J=\{0,1,\\2,\ldots,J, J+1\}$, and let $\states_K = \{0,1,2,\ldots,K, K+1\}$ represent the aggregated superstate space where $K\leq J$. If the state space is ordered, we define a function $s:S_J\mapsto S_K$ to be a valid state aggregation function if $s$ satisfies the following:
\begin{enumerate}
    \item If $i,j \in \states_J$ such that $i < j$, then $s(i) \leq s(j)$ for $s(i), s(j) \in \states_K$.
    \item If $i, j \in \states_J$ such that $j = i+1$, then, either $s(i) = s(j)$ or $s(j) = s(i) + 1$.
\end{enumerate}
\end{definition}
We also introduce the set $L_k = \{i \in \states_J: s(i) = k\}$ which is the set of states comprising superstate $k \in \states_K$.
For a policy based on an MDP with state space $\states_K$, $\pi_K$, we ``unaggregate'' the policy by stating that ${\pi}_J(\state) = {\pi}_K(k), \; \forall \state \in L_k, \; \forall k \in \states_K$.

\subsection{Aggregating transition probability matrices and rewards}

Suppose we have some TPM $P$ on state space $\states_J$, and let $s$ be some state aggregation function mapping $\states_J$ to $\states_K$ where $K, \, J$ are positive integers and $K \leq J$. We define $Q_K$ to be our aggregated TPM on state space $\states_K$ according to \cite{regnier2013state}. For any states $k, k' \in \states_K$ and action $a \in \actions$,
\vspace{-0.2cm}
\begin{equation}\label{eq. p to q}
    q_{kk'}^a = \frac{\sum_{\state \in L_k}\sum_{\state'\in L_{k'}}\beta_h p_{\state \state'}^a}{\sum_{\state \in L_k}\beta_h}. 
\end{equation}
where $\beta$ is a modified stationary distribution for $P$ where any absorbing state $h$ is modified such that $p_{hh}^a = 0$ and $p_{hh'}^a = \frac{1}{|S_J|-1}$ for all $\state' \neq \state$. See \cite{regnier2013state} and Appendix \ref{appendix: parameters} for details. In general, for most \glspl{tpm}, state aggregation leads to the loss of the Markov property \citep{kemeny1969finite}; however, this does not preclude the utility of an \gls{mdp} as the best (e.g., most interpretable, leading to the lowest regret, etc.) model for a given scenario. Because this paper is meant to consider the utility of an \gls{mdp} for an applied case study, we see the loss of the Markov property as being noteworthy, but not in direct contradiction to the goals of this paper.

Next, we can compute the aggregated rewards, $\rewards_Q$. In this study, we assume $r_Q(k,a)$ is equal to the simple average of $r(i, a)$ for every $i \in L_k$. That is, for each $a \in \actions$, $k \in \states_K$,
\vspace{-0.3cm}
\begin{equation}\label{eq. aggregate rewards}
    r_Q(k, a) = \frac{1}{|L_k|} \sum_{\state \in L_k}r(\state, a).
\end{equation} 
In some state aggregation frameworks \citep{bennouna2021learning}, it is necessary for $r_Q(k, a) = r(h_1, a) = r(h_2, a)$ for all $ h_1, h_2 \in L_k$. In our framework, we do not require this condition, which allows us to consider different types of problems where state aggregation can be useful.

\subsection{Estimation of transition probability matrices from data}
The TPM in an MDP can be estimated from observational data using maximum likelihood estimation (MLE) \citep{craig2002estimation}. Given observations from the MDP on the full state space $\states_J$, one can construct a $|\states|\times |\states| \times |\actions|$ observed count matrix $\obs$ where element $n_{\state \state'}^a$ is the number of recorded transitions from state $\state$ to state $\state'$ under action $a$ for each $\state, \state' \in \states$. The MLE of the TPM $P$ on the full state space is given by $\hat{P}$ with entries:

\begin{equation}\label{eq. Phat}
    \hat{p}_{hh'}^a = \frac{n_{hh'}^a}{\sum_{j = 0}^J n_{hj}^a},\quad \forall h, h'\in\states_J.
\end{equation}

A similar process is used to derive the aggregated model $\hat{Q}$, which is defined on the aggregated state space $\states_K$ for $K \leq J$ \cite{regnier2013state}. 
We can compute the probability of transitioning from state $k \in \states_K$ to state $k' \in \states_K$ under action $a$ using
\begin{equation}\label{eq. Qhat}
    \hat{q}_{kk'}^a = \frac{\sum_{j \in L_{k'}}\sum_{i \in L_k}  n_{ij}^a}{\sum_{j = 0}^J \sum_{i \in L_k} n_{ij}^a}.
\end{equation}

\subsection{The effects of aggregation policies in optimal stopping time MDPs}
Due to the conservation of state-ordering in our state aggregation function $s$, we can observe similarities between the properties of $\transprobs^a$ and $Q^a$.
\vspace{0.2cm}
\begin{proposition}\label{prop: Q is DFR}
    If $s:\states_J\mapsto \states_K$ is a valid state aggregation function, then TPM $P^a$ having the DFR property implies that the aggregated TPM $Q^a$ also has the DFR property.
\end{proposition}
\vspace{0.2cm}
\noindent We will defer this proof and all others to Appendix \ref{appendix: proofs}. Proposition \ref{prop: Q is DFR} will become useful when we consider how the structure of unaggregated optimal policies  relates to the structure of aggregated optimal policies. 

Consider an \gls{mdp} with action space $\actions = \{0, 1\}$ where 0 is the ``do-nothing" action and 1 is the ``intervene" action. There are sufficient conditions that guarantee the existence of an optimal policy that is a threshold policy \citep{alagoz2004optimal}. If we have an unaggregated MDP that meets these criteria, we can guarantee the following.
\vspace{0.2cm}
\begin{proposition}\label{prop: Q gives threshold policy}
Given an unaggregated MDP $(\states_J, \actions, P, \rewards)$ that satisfies the necessary conditions given by \cite{alagoz2004optimal} to guarantee a threshold policy and a state aggregation function which generates superstates $\{1,\ldots,K\}$ such that for all $k = 2, 3, \ldots, K-1, K$, we have that $|L_{k-1}| = |L_{k}|$ and 
$\frac{r_Q(k, \, 1) - r_Q(k-1, \, 1)}{r_Q(k \,, 1)} \leq \frac{r(\min(L_{k}), \, 1) -r(\max(L_{k-1}), \, 1)}{r(\min(L_{k}), \, 1)}$, the corresponding aggregated MDP $(\states_K, A, Q_K, R_{Q_K})$ is guaranteed a threshold policy that is optimal.
\end{proposition}
\vspace{0.2cm}
\noindent Thus, under certain conditions, we show that if the unaggregated model has a threshold policy, we also expect aggregated models to have a threshold policy. Note that the converse does not hold in general.

\begin{figure}[b]
    \centering
    \includegraphics[width=.9\textwidth]{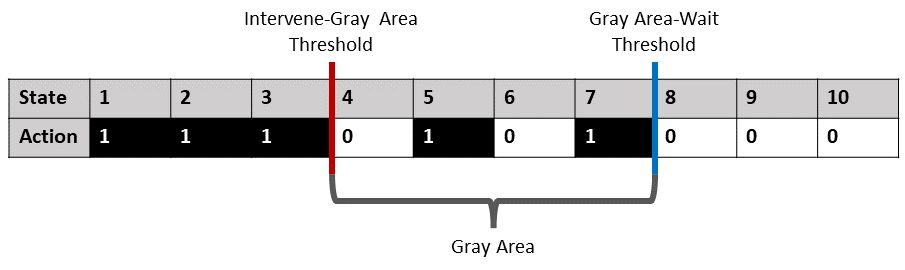}\vspace{-0.4cm}
    \caption{An illustration of a policy with a gray area where $\grayarea = \{4,5,6,7\}$. A \gls{dm} would likely expect the true threshold to exist within this set.}\label{fig: threshold/gray area}
\end{figure}
We now compare the ground truth model to the model estimated from data. Even if a ``ground truth" \gls{mdp} given by $(\states_J, \actions, P, \rewards)$ is guaranteed a threshold policy, its estimated \gls{mdp} given by $(\states_J, \actions, \hat{P}, \rewards)$ is \textit{not} guaranteed to have a threshold policy due to statistical error causing the sufficient conditions for a threshold policy to not hold. Solving a model with statistical error may result in a policy that has what we refer to as a \textit{gray area} or a region of the state space where the optimal policy alternates between actions $0$ and $1$ (potentially around the true threshold value). This gray area may be of particular concern to \glspl{dm} because the policy in this region is not intuitive. For example, if a \gls{dm} wishes to implement a threshold policy, the gray area would introduce ambiguity around where the best threshold lies (see Figure \ref{fig: threshold/gray area}). The gray area is formally defined as follows:

\begin{definition}\label{def: gray area} \textit{(Gray area)}
Let $\hat{\pi}^*$ be the estimated optimal policy of an MDP with a threshold policy that is optimal. Assume that $\hat{\pi}^*$ has at least one entry equal to $0$ and at least one entry equal to $1$. Let 
$$\Psi = \Big\{\state \in \{1, 2, \ldots, J-1\} : \, \hat{\pi}^*(\state) = 1, \hat{\pi}^*(\state + 1) = 0 \Big\}.$$
The \textbf{intervene-gray area threshold, $\grayarea_1$} is given by 
$$\grayarea_1 = \min \bigg\{\min \Big\{\state \in \{1, \ldots, J\}: \, \hat{\pi}^*(\state) = 0 \Big\}, \min \Psi \bigg\}.$$ The \textbf{gray area-wait threshold, $\grayarea_0$}, is given by 
$$\grayarea_0 = \max \bigg\{\max \Big\{\state \in \{1, \ldots, J\}: \, \hat{\pi}^*(\state) = 1 \Big\}, \max \Psi + 1 \bigg \}.$$ The \textbf{gray area, $\grayarea$} is given by  $$\grayarea = \{\state \in \states_J: \, \grayarea_1 < \state < \grayarea_0\}.$$
If $\grayarea = \varnothing$, we say that $\hat{\pi}^*$ is a threshold policy with threshold $T = \grayarea_1$.
\end{definition}
\noindent This definition is useful when discussing the structure of estimated optimal policies and how the distribution of the gray area is impacted by state space size and quantity of available data to parameterize the \gls{mdp}.

In summary, we have established the following results:
\begin{itemize}
    \item A \gls{tpm} $\transprobs$ having the DFR property implies that aggregated \gls{tpm} $Q$ will also have the DFR property (Proposition \ref{prop: Q is DFR})
    \item There are sufficient conditions that prove that an aggregated \gls{tpm} $Q$ will have a threshold policy, assuming the unaggregated \gls{tpm} $P$ also has a threshold policy (Proposition \ref{prop: Q gives threshold policy}).
    \item An estimated \gls{mdp} may not have a threshold policy due to statistical errors in estimated entries of the \gls{tpm}, either $P$ or $Q$.
\end{itemize}
We have also defined a measure of policy ambiguity, given by the gray area. These results inform our computational study described in the next section.

\section{Simulation Study of State Aggregation}\label{section: methods}
In this section, we describe a simulation study used to computationally investigate the effects of state aggregation on the performance of \gls{mdp} models. 
We consider an \gls{mdp} which extends \cite{regnier2013state}'s  \gls{dtmc} framework for a chronic disease progression by including a ``wait'' (0) and ``intervention'' (1) action. We simulate observations from the ground truth \gls{mdp} to generate synthetic observation counts $n_{ij}$ and then investigate how state aggregation affects the resulting policies generated by the estimated MDPs.

\subsection{Properties of the MDPs in this study}
Here, we state the important properties of the MDP structure. For the sake of space, we defer other details about the specific parameters and state aggregation function used in this study to Appendix \ref{appendix: parameters}. First, \cite{regnier2013state} gives us the following:

\begin{proposition}\citep[Proposition 1]{regnier2013state}\label{prop: P has DFR property}
The TPM $\transprobs^0$ described in Appendix \ref{appendix: parameters} has the DFR property. 
\end{proposition}

\noindent Given Proposition \ref{prop: P has DFR property}, the state aggregation procedure (Appendix \ref{appendix: parameters}), and Proposition \ref{prop: Q gives threshold policy}, the following can be proven:

\begin{proposition}\label{prop: Q gives threshold policy in our example}
Every MDP $(\states_J, A, P, R)$ and $(\states_K, A, Q_K, R_{Q_K})$ used in this study (detailed in  Appendix \ref{appendix: parameters}) is guaranteed to have a threshold policy that is optimal. 
\end{proposition} 
\noindent Proving that $(\states_J, \actions, \transprobs, \rewards)$ has an optimal threshold policy can be shown algebraically using the conditions in \cite{alagoz2004optimal}. Consequently, Proposition \ref{prop: Q gives threshold policy} proves that each MDP $(\states_K, \actions, Q_K, \rewards_{Q_K})$ in our study must have a threshold policy that is optimal. See Appendix \ref{appendix: proofs} for the formal proof.

\subsection{Simulation Procedures}
We now describe the simulation procedure for the experiments.
Our goal is to be able to compare the policies and remaining lifetimes estimates from \textit{estimated} MDP models (those with TPM $\hat{P}$ or $\hat{Q}_K$)  to those of the \textit{ground truth} MDP model (with TPM $P$). There are four steps to the experiment: 
\begin{enumerate}
    \item Generating synthetic observational data from the ground truth TPM $P^0$,
    \item Estimating the TPM for the MDP model using the synthetic data for a given state space $\states_J$ or $\states_K$,
    \item Solving the MDP model to obtain the estimated optimal policy, and
    \item Evaluating the estimated optimal policy in the ground truth model.
\end{enumerate}

\subsubsection{Generating synthetic observational data on system progression}
\hfill \break
First, we generate synthetic observational data from the ground truth model by simulating a positive integer $M$ system trajectories through the states according to the ground truth \gls{tpm} $P^0$.  Let $\state_{m,t}$ represent the state of system $m$ in period $t$ for $m = 1, \ldots, M$ and $t \in \timeset$. We assume $\state_{m, 0} = J \quad \forall m = 1, \ldots, M$. Given the state of system $m$ at time $t$, $\state_{m,t}$, we determine the state at time $t+1$ using Monte Carlo simulation. We sample a random variable, representing the next state, whose outcome is $j$ with probability $p^0_{h_{m,t},j}$. This process is repeated until time $t$ such that $\state_{m,t} = 0$, and this trajectory is finished. Because $P^0$ has the DFR property, each system entity is guaranteed to progress toward state $0$ as $t \rightarrow \infty$.
In our experiments, we consider $M\in\{10, 25, 50, 100, 500, 1000\}$, and we set $J=100$. We replicate the process of simulating $M$ system trajectories $R$ times, inputting the observations from all $M$ system trajectories of replication $r$ into an observed count matrix $\obs^r$ for $r=1,\ldots,R$. We use $R=100$ in our study.  

\subsubsection{Estimating transition probabilities and building the MDP model}\label{section: our model}
\hfill \break
Next, we use the synthetic observational data described above to construct MDP models with various levels of state aggregation. First, we solve the MDPs described by ($\states_J, \actions, \hat{P}, \rewards$) and ($\states_K, \actions, \hat{Q}_K, \rewards_{Q_K}$) to obtain optimal policies $\hat{\pi}_{\hat{P}}^*$ and $\hat{\pi}_{\hat{Q}_K}^*$, respectively using the Python \verb|mdptoolbox| module \citep{chades2014mdptoolbox}. We denote the optimal policy of any arbitrary estimated model as $\hat{\pi}^*$ with no subscript.  After these policies are generated, we 
investigate the optimal policy structure, different definitions of ambiguity (e.g., gray area, lack of threshold precision due to state aggregation, etc.), and how modeling decisions impact policy ambiguity and different value metrics of the policy. 

First, we analyze the resulting policies to determine whether aggregation affects when intervention is recommended. By Proposition \ref{prop: Q gives threshold policy in our example}, the ground truth MDP described by ($\states_J, \actions, \transprobs, \rewards$) has an optimal threshold policy, and we compute that the optimal threshold is $T_P = 24$. Hence, if there were no statistical errors in the TPM, then 100\% of the estimated policies in the $R$ replications would recommend intervention for $h \leq T_P$ and waiting for $h>T_P$. 
 However, the absence of a sufficient number of observations when constructing an estimated TPM can lead to ambiguity in the threshold. This can be shown when different replications of an experiment for fixed $M$ and $K$ yield different optimal policies. 
 
 

Next, we consider policy structure and ambiguity. By Proposition \ref{prop: Q gives threshold policy in our example}, the MDP described by ($\states_J, \actions, \transprobs, \rewards$) has an optimal threshold policy. However, in general, the MDPs with TPMs $\hat{P}$ and $\hat{Q}$ are not guaranteed to have threshold policies that are optimal, due to statistical error. We investigate whether these policies have a threshold structure and, if not, consider the size and empirical distribution of the gray area over the 100 replications. 

Finally, we make value comparisons between the ground truth and the models estimated from data. One value metric is to quantify the expected losses in value incurred from estimation and/or aggregation, called \textit{expected regret}. Let $v(\hat{\pi}^*, \cdot)$ denote the value of policy $\hat{\pi}^*$ given by  \eqref{eq. value to go} where the second parameter is the TPM in which $\pi^*$ is evaluated, either $P, \hat{P}$, or $\hat{Q}$. We calculate expected regret $\xi$ using
 \begin{equation}\label{eq. expected regret}
    \xi(\hat{Q}_K) = \val(\pi_P^*, P) - \val(\hat{\pi}_{\hat{Q}_K}^* P).
 \end{equation} 

Some \glspl{dm} may prefer a state space $\states_K$ with TPM $\hat{Q}_K$ that minimizes $\xi$.

\subsection{Threshold policy assumption}
It is possible that a modeler would create an MDP under the assumption that the MDP has a threshold policy, and he or she would only consider policies with a threshold structure. As an addendum to the aforementioned experiment, we will also compute and evaluate estimated optimal policies from estimated TPMs under the assumption of a threshold policy. We will denote the optimal policy under the threshold assumption as $\bar{\pi}^{*}$. Let $\mathfrak{T}$ be the set of all threshold policies, and let $v$ be the value function. Then,
    $\bar{\pi}^{*} = \argmax_{\pi \in \mathfrak{T}} \, v(\pi)$, which can be solved by iterating over all possible $\bar{\pi} \in \mathfrak{T}$. Note that $\pi_P^* = \bar{\pi}_P^{*}$, since the true optimal policy is a threshold policy by Proposition \ref{prop: Q gives threshold policy in our example}.

\section{Results}\label{section: results} 

Here, we discuss the interpretability, correctness, and value of the estimated optimal policies computed from differently aggregated \glspl{mdp} with \glspl{tpm} estimated from data simulated from a single ground truth \gls{tpm} (Appendix \ref{appendix: parameters}).



\subsection{Implications of state aggregation on intervention recommendations}

\begin{figure}[ht]
\centering
\begin{subfigure}{.32\textwidth}
  \centering
\includegraphics[width=.98\linewidth]{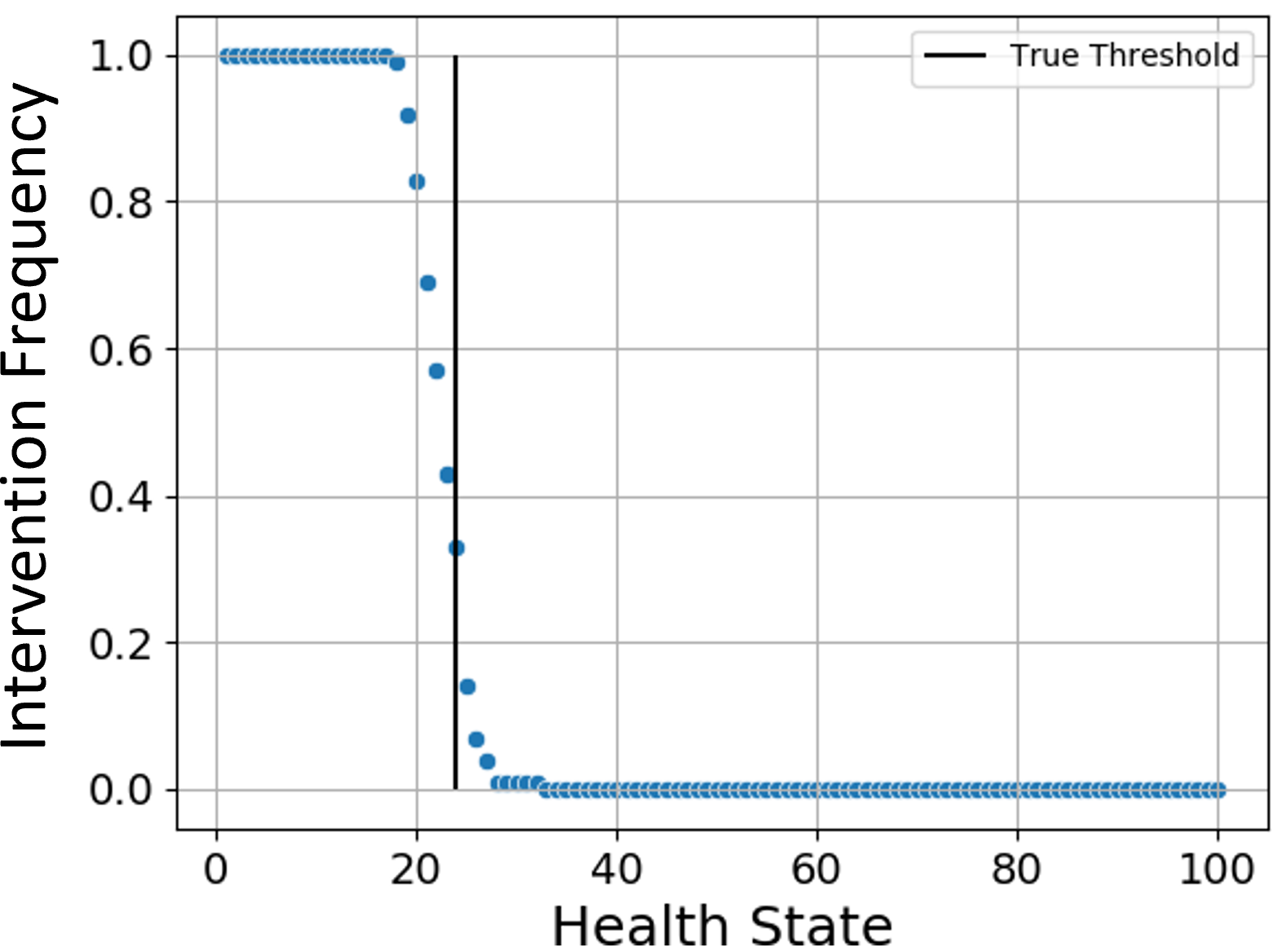}  
  \caption{$K = 100$, $M = 50$ (threshold)}
  \label{fig: subfig: policies - thresh, no lumping, 50 traj}
\end{subfigure}
\begin{subfigure}{.32\textwidth}
  \centering
    \includegraphics[width=.98\linewidth]{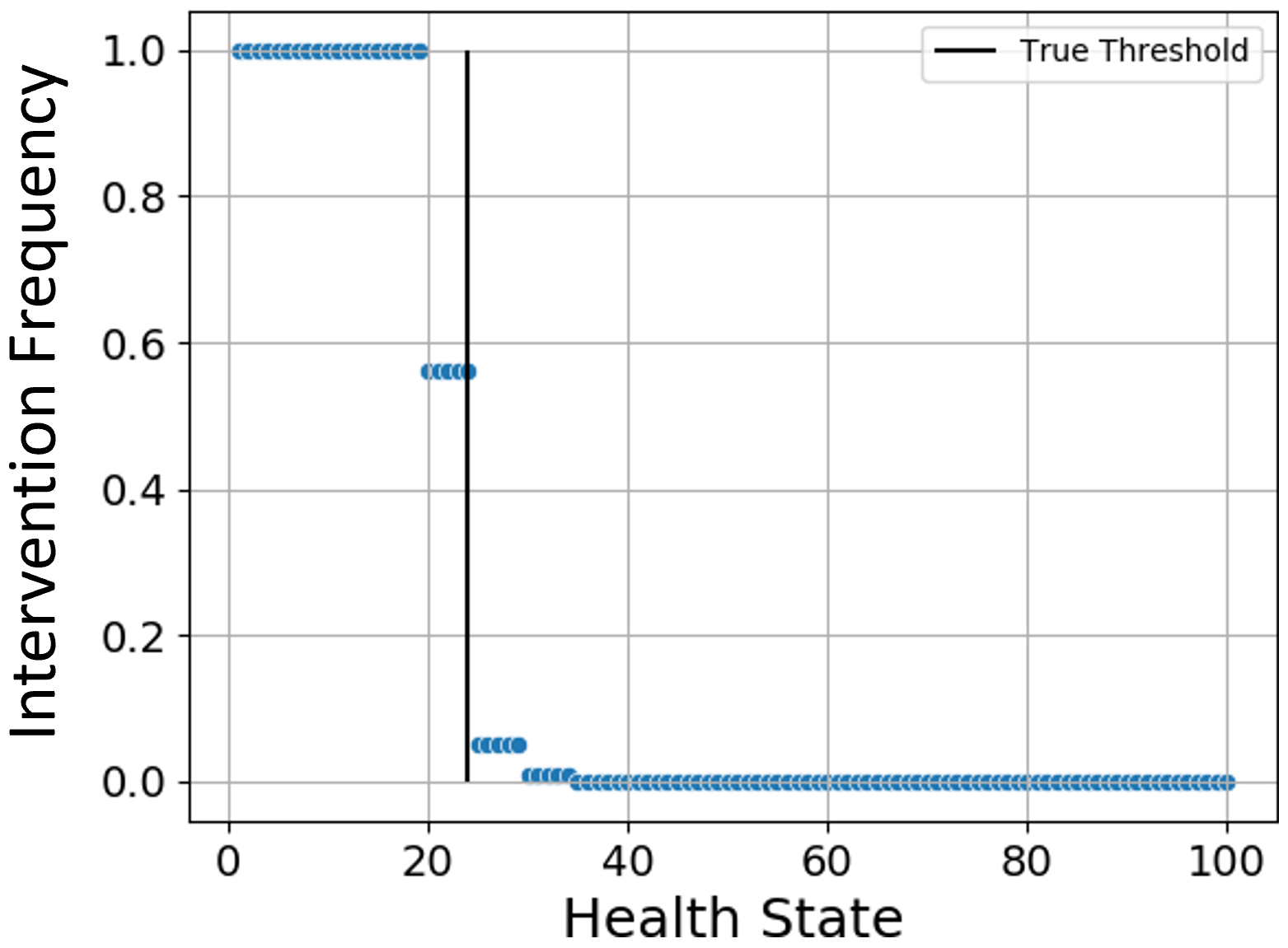}  
  \caption{$K = 20$, $M = 50$ (threshold)}
  \label{fig: subfig: policies - thresh, 20 states, 50 traj}
\end{subfigure}
\begin{subfigure}{.32\textwidth}
  \centering
    \includegraphics[width=.98\linewidth]{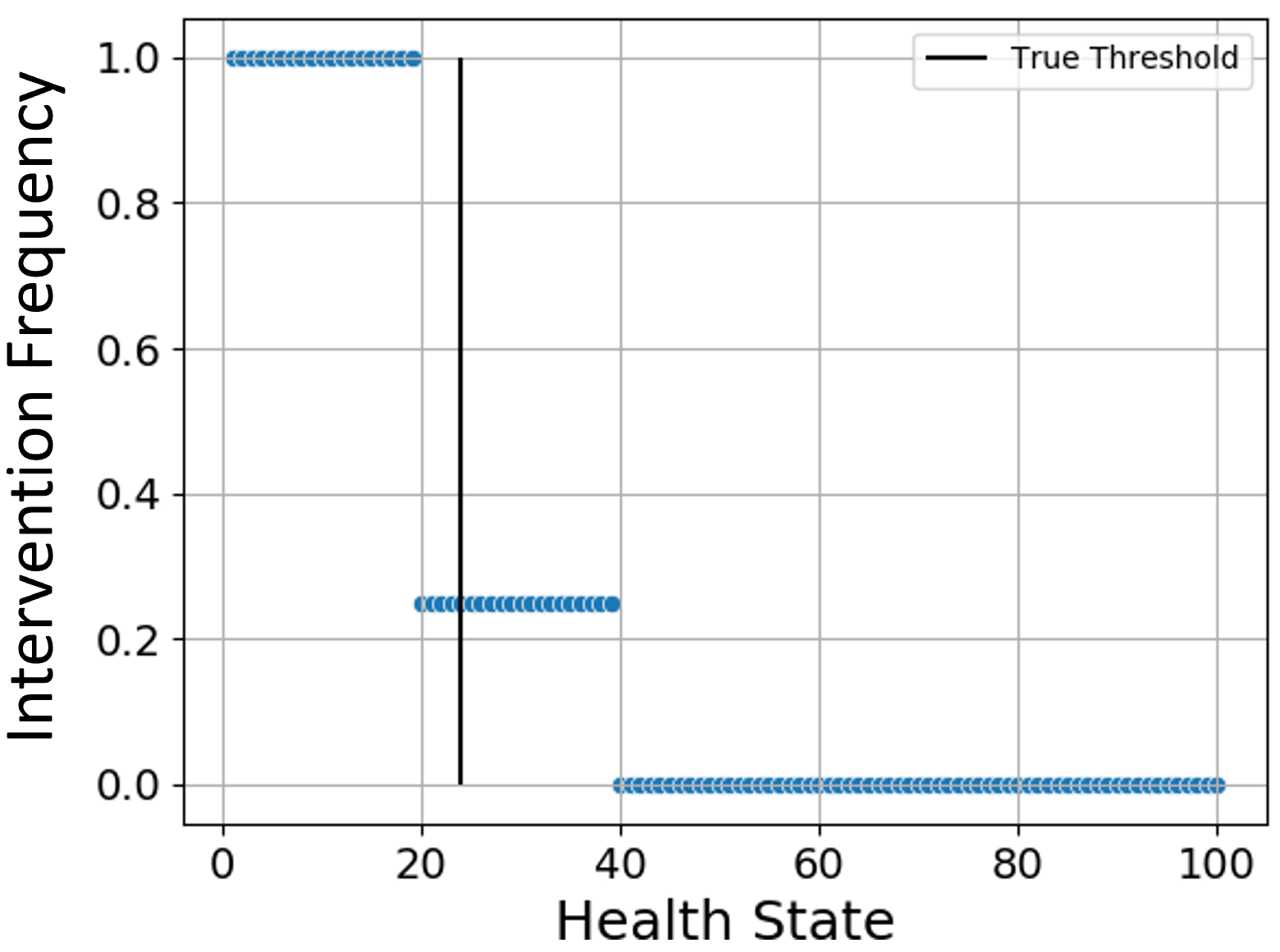}  
  \caption{$K = 5$, $M = 50$ (threshold)}
  \label{fig: subfig: policies - thresh, 5 states, 50 traj}
\end{subfigure}
\\
\begin{subfigure}{.32\textwidth}
  \centering
\includegraphics[width=.98\linewidth]{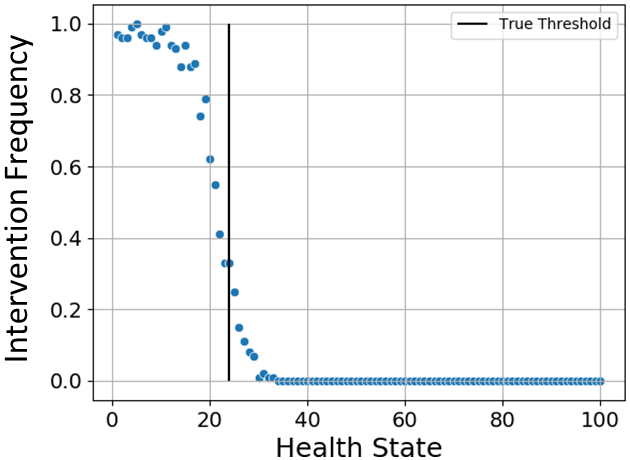}
  \caption{$K = 100$, $M = 50$}
  \label{fig: subfig: policies - no lumping, 50 traj}
\end{subfigure}
\begin{subfigure}{.32\textwidth}
  \centering
    \includegraphics[width=.98\linewidth]{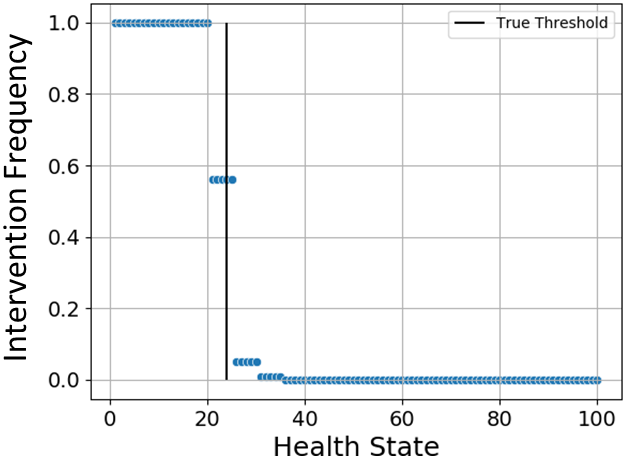}
  \caption{$K = 20$, $M = 50$}
  \label{fig: subfig: policies - 20 states, 50 traj}
\end{subfigure}
\begin{subfigure}{.32\textwidth}
  \centering
    \includegraphics[width=.98\linewidth]{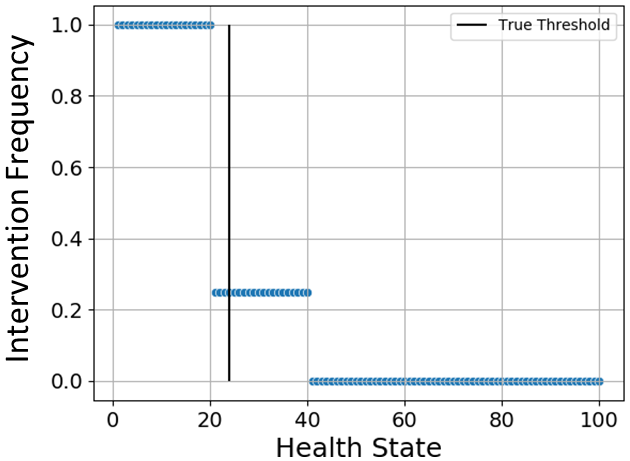}  
  \caption{$K = 5$, $M = 50$}
  \label{fig: subfig: policies - 5 states, 50 traj}
\end{subfigure}
\\
\begin{subfigure}{.32\textwidth}
  \centering
\includegraphics[width=.98\linewidth]{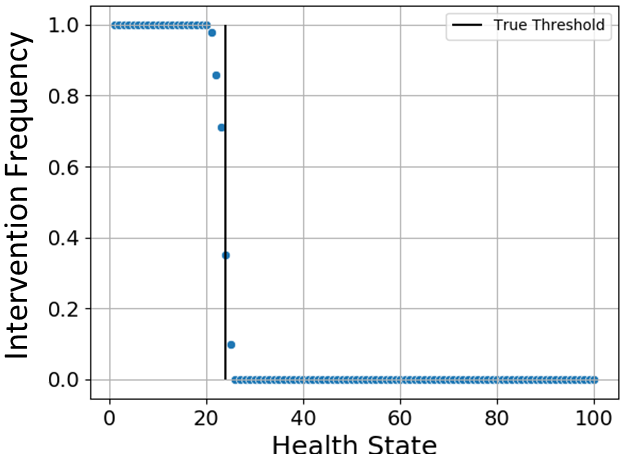}
  \caption{$K = 100$, $M = 500$}
  \label{fig: subfig: policies - no lumping, 500 traj}
\end{subfigure}
\begin{subfigure}{.32\textwidth}
  \centering
    \includegraphics[width=.98\linewidth]{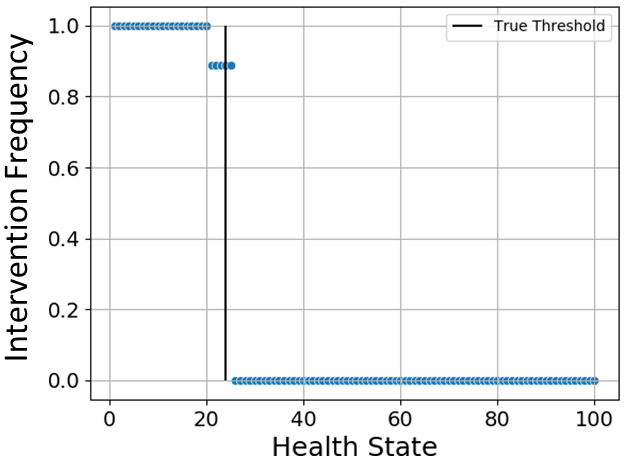}
  \caption{$K = 20$, $M = 500$}
  \label{fig: subfig: policies - 20 states, 500 traj}
\end{subfigure}
\begin{subfigure}{.32\textwidth}
  \centering
    \includegraphics[width=.98\linewidth]{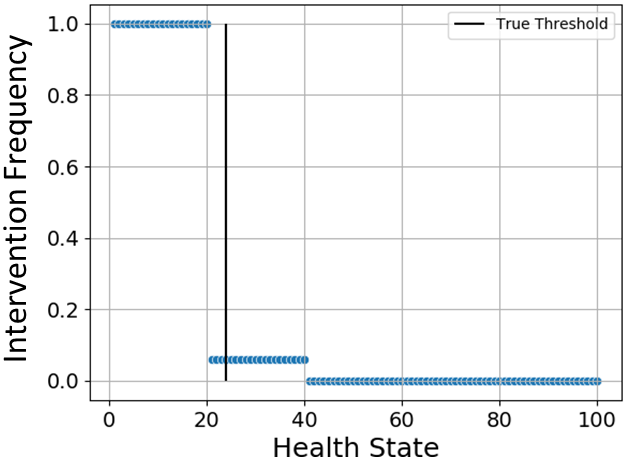}  
  \caption{$K = 5$, $M = 500$}
  \label{fig: subfig: policies - 5 states, 500 traj}
\end{subfigure}
\caption{The frequency of the recommendation to ``Intervene” in each state across 100 replications.   Each column corresponds to a state space with $K=100$, $K=20$, and $K=5$, respectively with 50 available system trajectories (top two rows) and 500 available system trajectories (bottom row). A threshold policy is assumed in Figures 3(a) - 3(c). The optimal threshold is shown by the vertical line at state 24. 
\label{fig: sampled policy recommendations}}
\end{figure}

We first consider how state aggregation and data availability impact the states in which intervention is recommended and at what frequency. Figure \ref{fig: sampled policy recommendations} shows the observed frequency of the intervention recommendation in each state as the availability of data and level of aggregation is varied. 

First, we observe the differences in the frequency of the intervention recommendation  in each state across $R$ replications when there is no threshold assumption (Figures \ref{fig: subfig: policies - no lumping, 50 traj} - \ref{fig: subfig: policies - 5 states, 500 traj}). We observe that as the state space size decreases, the intervention frequency plots approach an identical threshold policy for every replication.
A similar trend occurs when we increase the number of available system trajectories, $M$. As $M$ increases, these $R$ estimated policies begin to converge to a single threshold policy, although not necessarily the same threshold policy recommended by the true optimal policy $\pi_P^*$.
For example, take $K=5$ in Figures \ref{fig: subfig: policies - 5 states, 50 traj} and \ref{fig: subfig: policies - 5 states, 500 traj}. As the number of observed trajectories $M$ increases from $M=50$ to $M=500$, the intervention recommendations approach a threshold $\hat{T}_{\hat{Q}_5} = 20$. 

While using a highly aggregated model may be useful for very low $M$ values, we see here that the lack of state precision can lead to early or late intervention. We compare this to Figures \ref{fig: subfig: policies - no lumping, 50 traj} and \ref{fig: subfig: policies - no lumping, 500 traj} which show $K=J=100$. As $M$ increases from $50$ to $500$, the region of states for which there is disagreement about the optimal action shrinks from states $1, \ldots, 36$ to states $21, \ldots, 25$. We expect that as $M$ approaches infinity, the plot would converge toward a threshold policy plot with threshold $T_P$. Hence, there is a high potential for policy improvement under $\states_J$ as available data increases, whereas there is a limit to how much a policy under $\states_K$ can improve due to the lack of state precision.

Now, we observe how the assumption of a threshold policy changes the variability in the 100 estimated optimal policies. Consider the top two rows of Figure \ref{fig: sampled policy recommendations}, where the top row considers $M = 50$ system trajectories using a threshold policy assumption and the second row considers $M = 50$ system trajectories without the threshold assumption. We can observe that the threshold assumption leads to much better estimation of the true optimal policy in the $K = 100$ case, although there is no easily observable difference in the intervention recommendations for the $K = 20$ or $K = 5$ case. In Figure \ref{fig: subfig: policies - thresh, no lumping, 50 traj}, the region of states for which there is disagreement about the optimal action is states $18, \ldots, 33$, compared Figure \ref{fig: subfig: policies - no lumping, 50 traj} in which the ambiguous region consists of states $1, \ldots, 33$. For unaggregated state spaces, the threshold assumption appears to mitigate the impacts of statistical uncertainty in the lower states, which tend to have fewer observations.

\subsection{Implications of state aggregation on the distribution of the gray area}

\begin{figure}[h]
\centering
\includegraphics*[width=.6\textwidth,height= .6\textheight,keepaspectratio]{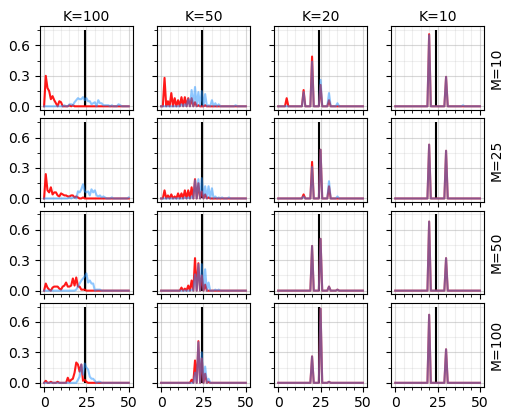}
\caption{Distributions of the intervene-gray area thresholds (red) and the gray area-wait thresholds (blue) for several pairings of number of system trajectories, $M$, and number of states, $K$.\label{fig: thresholds treat and wait}}
\end{figure}

Now, we consider the effect of state aggregation on the empirical distributions of the intervene-gray area and gray area-wait thresholds in the estimated optimal policies obtained from 100 replications of the experiment.
In Figure \ref{fig: thresholds treat and wait}, we plot the frequency of the locations of $\grayarea_0$ (the critical state above which $\hat{\pi}^*(\state) = 0$) and $\grayarea_1$ (the critical state below which $\hat{\pi}^*(\state) = 1$). Note that the state space is truncated at state 50 since there is no ambiguity in the optimal policy after that point. Furthermore, the ``spikes" that we observe for the subplots where $K < 100$ are due to the fact that a threshold can only exist at specific states for aggregated models (e.g., for K=10, our threshold can only exist at states that are multiples of 10).

The gray area in Figure \ref{fig: thresholds treat and wait} can be thought of as the space between the red and blue distributions representing $\grayarea_1$ and $\grayarea_0$, respectively. From this figure, we see that as the number of system trajectories increases, the distributions of $\grayarea_1$ and $\grayarea_0$ approach each other (i.e., the gray area shrinks). The same convergence occurs as the number of states in the state space decreases. For $K=20$ and $K=10$, we observe a very tight overlap in the distribution of the two thresholds, which lends the assumption that many of the 100 replications yield a threshold policy for these state spaces. However, the leftward bias we observe for low values of $M$ likely indicates that the computed threshold, $\hat{T}_{\hat{Q}}$ is often less than the true optimal threshold, $T_P$.

Another important observation from Figure \ref{fig: thresholds treat and wait} is how quickly the two empirical threshold distributions obtained from the 100 replications converge. We see that for $K=50$, the two distributions shift from almost complete separation to a high area of intersection between $M=10$ and $M=25$. This could be explained by the finding from  \cite{regnier2013state} in showing that statistical error decreases quickly even with small increases in $M$.
Hence, the optimal policies using high values of $K$ (i.e., a larger, less aggregated state space) will show quicker reductions in the ambiguous gray area than those with low values of $K$ (i.e., a smaller, more aggregated state space) with only small increases in $M$. 

\subsection{Implications of state aggregation on remaining system lifetime estimates}\label{section: remlife}

\begin{figure}[ht]
\centering
\begin{subfigure}{.48\textwidth}
  \centering
    \includegraphics[width=.95\linewidth]{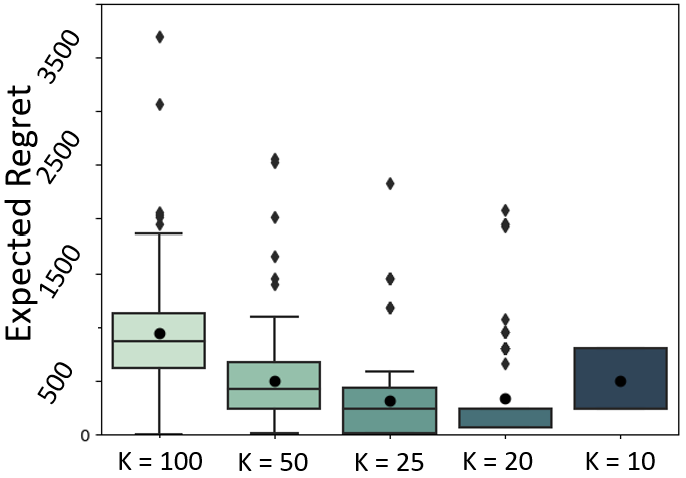}  
  \caption{25 system trajectories observed}
  \label{fig: subfig: regret M=25}
\end{subfigure}
\begin{subfigure}{.48\textwidth}
  \centering
  \includegraphics[width=.95\linewidth]{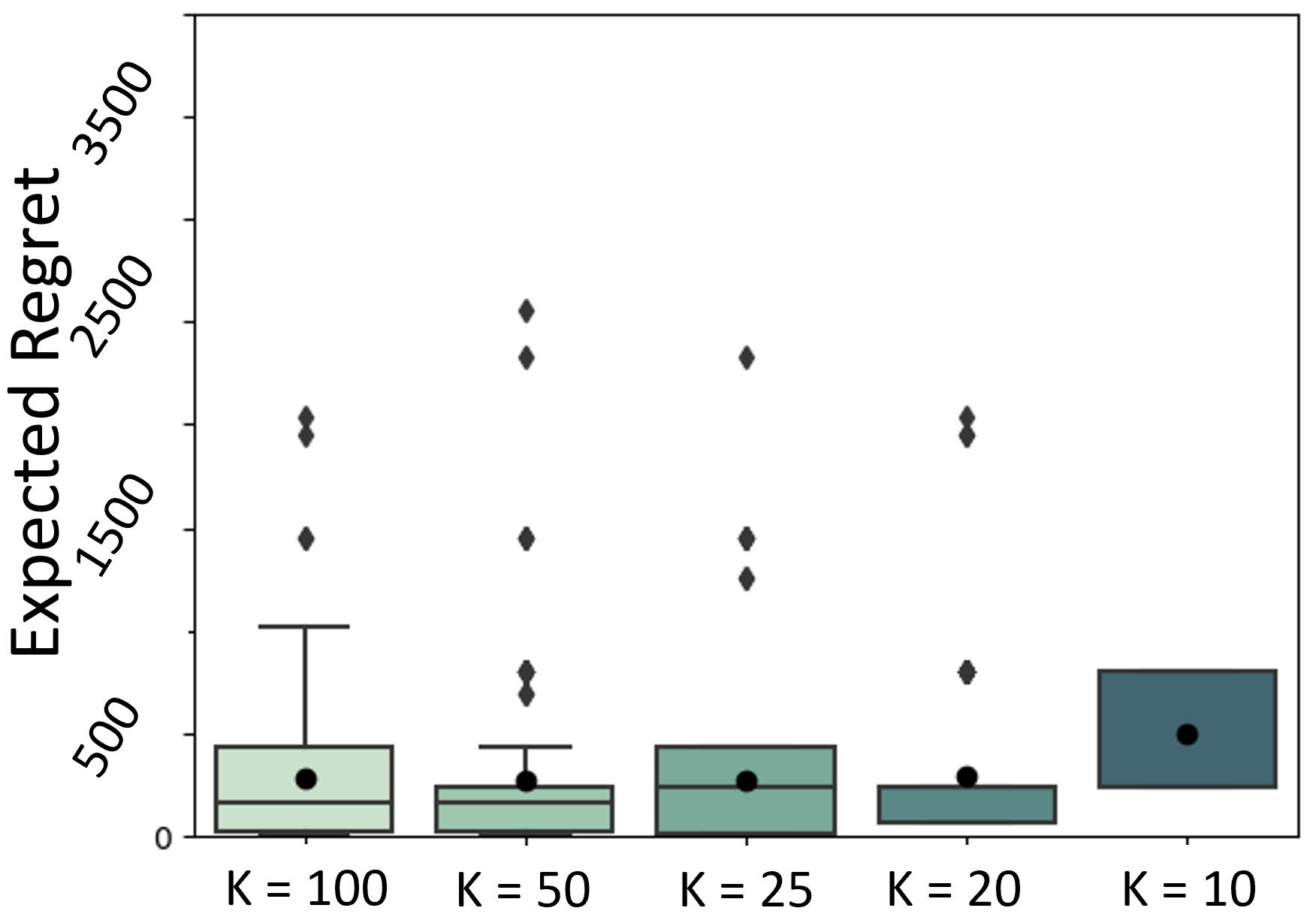}  
  \caption{25 system trajectories observed (threshold)}
  \label{fig: subfig: regret M=25 - thresh}
\end{subfigure}
\\
\begin{subfigure}{.48\textwidth}
  \centering
  \includegraphics[width=.95\linewidth]{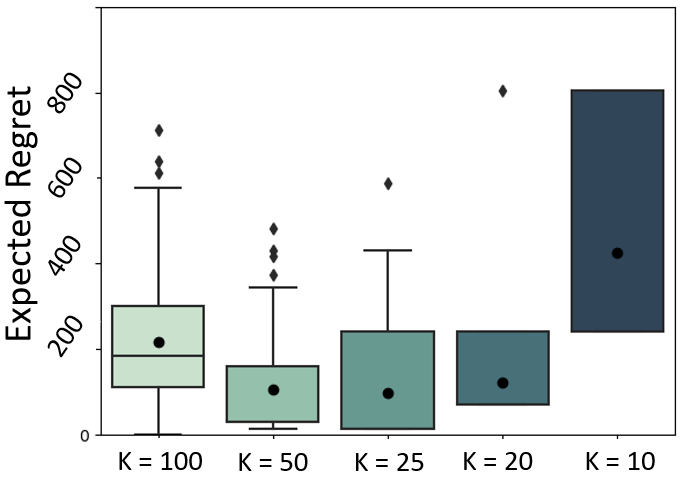}  
  \caption{100 system trajectories observed}
  \label{fig: subfig: regret M=100}
\end{subfigure}
\begin{subfigure}{.48\textwidth}
  \centering
  \includegraphics[width=.95\linewidth]{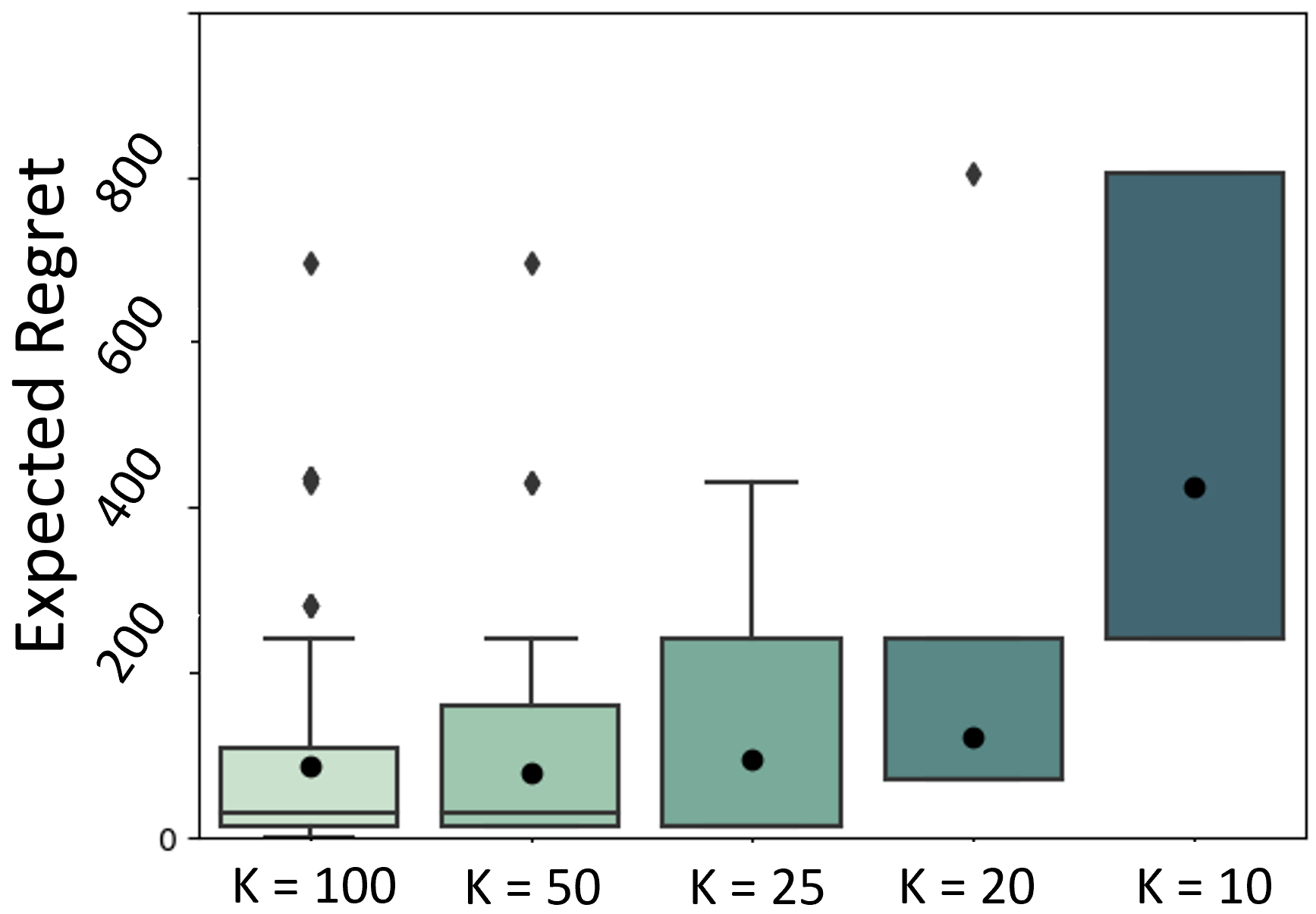}  
  \caption{100 system trajectories observed (threshold)}
  \label{fig: subfig: regret M=100 - thresh}
\end{subfigure}
\caption{The expected regret from making intervention decisions from the optimal policy of an estimated model as a function of the number of system trajectories, $M$. The sample mean of the expected regret for each model is represented by a dot. \label{fig: expected regret}}
\end{figure}

In this subsection, we investigate the impacts of state aggregation and data availability on the expected regret, defined in \eqref{eq. expected regret}. Figure \ref{fig: expected regret} shows how the empirical distribution of expected regret across the 100 replications changes as the number of system trajectories increases. The sample mean of the expected regret for each model is represented by a dot.  

First, consider Figures \ref{fig: subfig: regret M=25} and \ref{fig: subfig: regret M=100}. We note a rough ``U"-shaped curve in the sample mean of the expected regret as the number of aggregated states decreases. For experiments without a threshold policy assumption, we observe some moderate level of state aggregation $K$ which minimizes the mean observed expected regret for a given  number of system trajectories available, $M$.
 
 Figures \ref{fig: subfig: regret M=25} and \ref{fig: subfig: regret M=100} show how the empirical distribution of expected regret across the 100 replications changes as the number of system trajectories increases. The sample mean of the expected regret for each model is represented by a dot.  In both Figures \ref{fig: subfig: regret M=25} and \ref{fig: subfig: regret M=100}, we note a rough ``U"-shaped curve in the sample mean of the expected regret as the number of aggregated states decreases. In each case, we observe some moderate level of state aggregation (in these cases, both $K = 25$) that minimizes the mean observed expected regret for a given  number of system trajectories available, $M$.

As we might expect, as the number of observed trajectories $M$ increases, the expected regret for a given model with fixed $K$ either decreases or remains constant due to more accurate estimates of the TPM. However, the rate at which the expected regret for each model decreases is not the same for each level of aggregation. From Figure \ref{fig: subfig: regret M=25} to Figure \ref{fig: subfig: regret M=100}, we see that the expected regret for the models defined by $K=10$ remain relatively unchanged as $M$ increases (note that y-axis scales differ on each row); the mean expected regret decreases from approximately 505  months per thousand system trajectories to approximately 426 months per thousand system trajectories, about a 15.6\% decrease.
On the other hand, $\hat{P}$, defined by $K=100$, decreases comparatively rapidly as $M$ increases. The mean expected regret for $K = 100$ decreases from approximately 948 months per thousand system trajectories to approximately 218 life months per thousand system trajectories, which is about a 77\% decrease. For reference, the ``never intervene" policy yields an expected regret of over $25,000$ months per thousand system entities, which shows that intervention vastly improves the life of the system entity in this scenario and demonstrates the scale of improvements from moderate aggregation.

Based on Figures \ref{fig: subfig: regret M=25} and \ref{fig: subfig: regret M=100}, one might conjecture that the lack of precision in highly aggregated models (e.g., $\hat{Q}_{10}$) will lead to underperformance relative to more finely aggregated models (e.g., $\hat{P}$), especially when the number of observed trajectories, M, is large. However, in some cases, the impact of the lack of precision depends on the value of $T_P$ and its location relative to the upper limit cutoff states in $\states_J$ for each aggregated state $h \in \states_K$. For example, if $J=100$ and $K=5$, the cutoff states would be $\{20, 40, 60, 80, 100\}$. 
Suppose our model had $T_P = 20$ instead of $T_P = 24$. Because $20$ is a cutoff state for $\states_5$, there will be no early or late intervention when a threshold policy is obtained for $K=5$. Hence, our model with $K=5$ will likely attain $0$ regret by a relatively low $M$ value, in which case $\hat{P}$ can, at best, have an equal expected regret. Deviations from the aforementioned U-curve with regard to the sample mean of the expected regret can also be partially attributed to the relationship between $T_P$ and state cutoffs. See Appendix \ref{appendix: alternative TPM} for an example of this.

Next, we consider how the addition of a threshold policy assumption alters the expected regret. The most noticeable difference between Figures \ref{fig: subfig: regret M=25} and \ref{fig: subfig: regret M=100} and Figures \ref{fig: subfig: regret M=25 - thresh} and \ref{fig: subfig: regret M=100 - thresh} is that the expected regret for $K = 100$ and $K =50$ is substantially lower when we assume a threshold policy. In \ref{fig: subfig: regret M=25 - thresh}, the mean expected regret for $K = 100$ is approximately 285 months, whereas in \ref{fig: subfig: regret M=100 - thresh}, the mean expected regret is approximately 88 months. Both values are substantially lower than the mean expected regrets observed in the case where a threshold policy is not assumed (Figures \ref{fig: subfig: regret M=25} and \ref{fig: subfig: regret M=100}), which have expected regrets of approximately 948 and 218 months, respectively. This result is intuitive for two reasons. First, the true optimal policy is also a threshold policy, so this assumption does not lead to an inherent disadvantage when modeling. Second, the main drawback of large state spaces is the relatively low number of system observations to draw from in each state. By forcing the estimated optimal policy to have a threshold structure, we are not simply choosing which action to take in which state, we are choosing the threshold location. This suggests that threshold policies tend to be robust against errors in entries $\hat{p}_{hh'}$ in $\hat{P}$. Hence, state aggregation loses its advantage when we already employ the threshold policy assumption. For highly aggregated models, assuming a threshold policy has a lesser effect on decreasing regret. In fact,  we achieve near-identical results for $K = 10$ in Figures \ref{fig: subfig: regret M=25} and \ref{fig: subfig: regret M=25 - thresh} and Figures \ref{fig: subfig: regret M=100} and \ref{fig: subfig: regret M=100 - thresh}.

\begin{figure}[ht]
\centering
\begin{subfigure}{.49\textwidth}
  \centering
    \includegraphics[width=.99\linewidth]{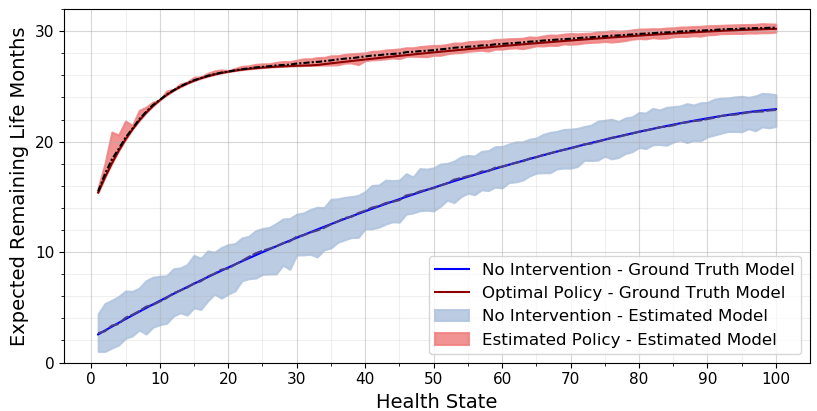}  
  \caption{100 states ($\hat{P}$), $M = 25$}
  \label{fig: subfig: rem lifetime K=100}
\end{subfigure}
\begin{subfigure}{.49\textwidth}
  \centering
  \includegraphics[width=.99\linewidth]{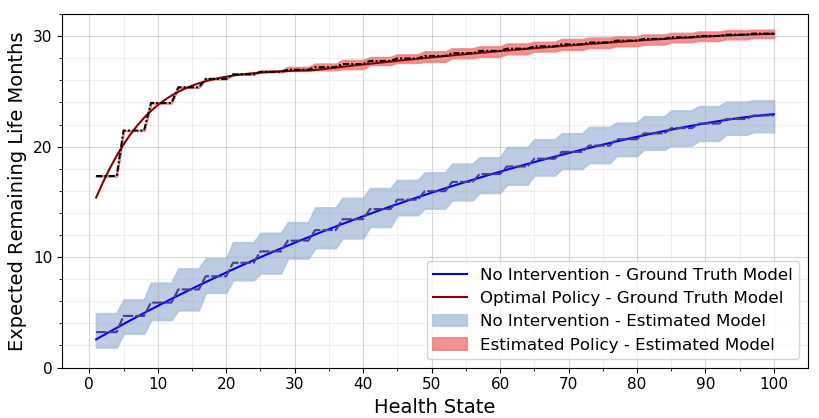}  
  \caption{25 aggregated states ($\hat{Q}_{25}$)}
  \label{fig: subfig: rem lifetime K=25}
\end{subfigure}
\\
\begin{subfigure}{.49\textwidth}
  \centering
  \includegraphics[width=.99\linewidth]{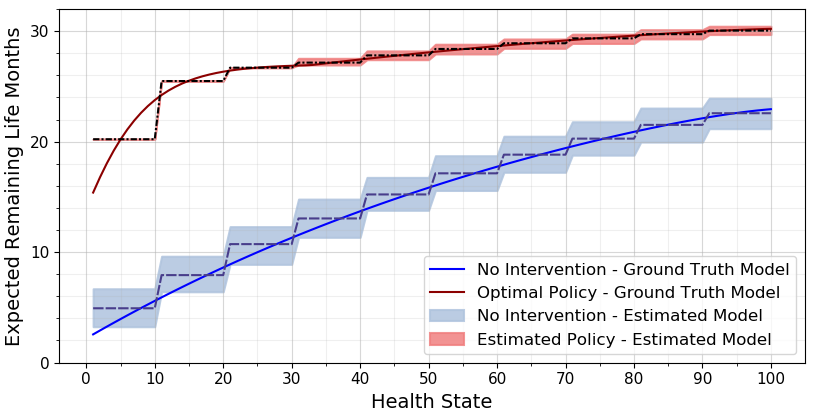}  
  \caption{10 aggregated states ($\hat{Q}_{10}$), $M = 25$}
  \label{fig: subfig: rem lifetime K=10}
\end{subfigure}
\begin{subfigure}{.49\textwidth}
  \centering
  \includegraphics[width=.99\linewidth]{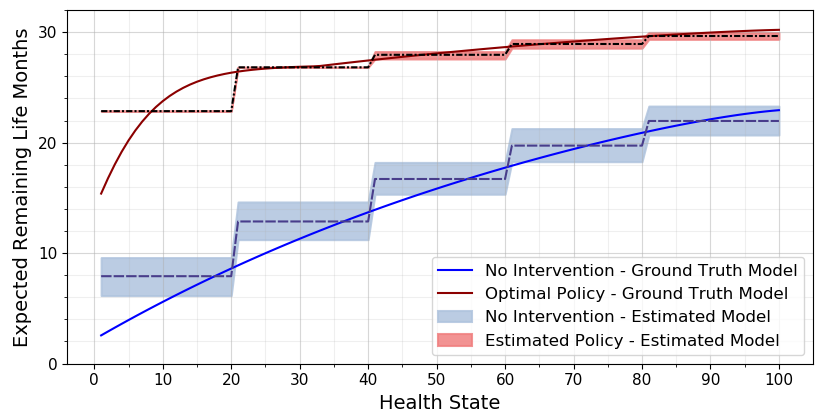}
  \caption{5 aggregated states ($\hat{Q}_5$), $M = 25$}
  \label{fig: subfig: rem lifetime K=5}
\end{subfigure}
\caption{The 95\% confidence interval for the estimated expected remaining system lifetime when the system begins in each health state under the estimated optimal policy (red) and no-intervention policy (blue). \label{fig: rem lifetime}}
\end{figure}

Lastly, we examine the difference between the expected remaining life months using the optimal policy $\hat{\pi}^*$ versus the no-intervention policy $\pi_0$. Figure \ref{fig: rem lifetime} shows the $95\%$ confidence interval for the estimated remaining life months for each state according to the \gls{dm} (i.e., without access to the ground truth), given by $\val(\pi_0, \hat{Q}_K)$ and $\val(\hat{\pi}_{\hat{Q}_K}^* \hat{Q}_K)$.

We first observe that 
the difference between the sample means of $v(\hat{\pi}_{\hat{Q}_K}^*, \hat{Q}_K)$ and $v(\pi_0, \hat{Q}_K)$ at $\state=100$ increases as $K$ decreases.
Secondly, we note some biases in the $\hat{Q}$ estimates. One such bias is that for aggregated policies in lower states, the no-intervention policy often overestimates the expected remaining lifetime. Furthermore, due to the sharper increase of the optimal value function for low states, the estimated optimal value function for aggregated models \ref{fig: subfig: rem lifetime K=10} and \ref{fig: subfig: rem lifetime K=5} will only closely approximate the true optimal value for the states near the midpoint of the aggregated state. Once the value function becomes more level at approximately state 20, the estimated expected remaining lifetime under the estimated optimal policy is a closer estimate of the true expected remaining lifetime, although estimated remaining lifetime in higher states are slightly underestimated.  

\section{Conclusion}\label{section: conclusion}
In this study, we 
show that estimated \gls{mdp} models with moderately aggregated state spaces can generate policies that lead to less decision ambiguity and lower expected regret than coarsely or finely aggregated models. This finding differs from the recommendation from \cite{regnier2013state} that \gls{dtmc} models of disease prognosis under no treatment should use little to no state aggregation. We find that, even for $1000$ observations, the mean expected regret for taking the optimal policy estimated from a model using the full 100-state state space (i.e., using $\hat{P}$) was still higher than those using aggregated TPMs $\hat{Q}_{25}$ and $\hat{Q}_{50}$. This finding reveals that minimizing error in the TPM does not necessarily translate to minimizing the expected regret of the MDP model. 

Furthermore, we find that for a low number of available system trajectories, moderately aggregated MDPs are more likely to closely approximate the true threshold policy and eliminate some of the ambiguity from the gray area. However, it is still possible to over-aggregate a state space, where precision losses lead to poor estimation of the threshold. Overall, one of the key problems with choosing the most desirable state space to model the system's progression with a threshold policy is the ``luck" of where that true unknown threshold $T_P$ lies on the aggregated state space. 
Hence, the generalization that the full model's estimated \gls{tpm} will lead to lower regret than an aggregated model's estimated \gls{tpm} as the number of observed system trajectories increases is not necessarily true. Future work might consider how to estimate this threshold location from limited data.

Our finding that moderately aggregated \glspl{tpm} can perform better with regard to threshold estimation and expected regret does not hold when the estimated optimal policy is required to be a threshold policy. When this assumption is made, $\hat{P}$ often leads to a lower mean expected regret in our computational study than $\hat{Q}_K$ for $K < 100$. Because the true optimal policy is a threshold policy, we incur no disadvantage by making this assumption. Furthermore, requiring a threshold policy protects against individual outlier point estimates for $P^0$, particularly for low values of $h$. The natural question is whether or not it makes sense to assume a threshold policy in a general setting where we are unsure of the true optimal policy structure. The answer to this question would depend on the system being modeled and the preferences of the \glspl{dm}. For example, if the ``system" was a patient being treated for a chronic disease, threshold policies are commonly assumed. From a clinician perspective, threshold policies are intuitive and easier to implement than non-threshold policies. Hence, a threshold policy assumption may be preferable, even if the underlying disease progression model would not lead to a threshold policy that is optimal. For other systems, a threshold policy may not be an advantageous assumption.

Our work is not without limitations. First, we consider a hypothetical deteriorating Markov process from literature \citep{regnier2013state}. The use of hypothetical models has been used before to enable comparisons to the ``ground truth" optimal policy and remaining lifetime estimates that would not be available otherwise \citep{regnier2013state,moreno2010optimal}. 
However, the model used was one of many possible models that could be considered with the desired properties specified in \S\ref{subsection:mdps of interest}. Future work could consider other TPM designs using the framework from \cite{regnier2013state} or consider different structures of TPMs altogether, such as models in which the system can never improve. These different structures may yield different results in terms of the utility of state aggregation. 

Our analysis motivates opportunities for future research. Future work could investigate methods for comparing different state space designs in the absence of a ground truth model. Our work suggests that the choice of state space can indeed influence the quality of the resulting recommendations and remaining system life estimates. A rigorous approach to comparing different potential state space designs when there is no ground truth model could be beneficial to modelers. Future work could also consider the influence of censored observational data (e.g., sparsity in transitions between high-health and low-health states in chronic disease models) on state space design. 
Furthermore, this article focused on optimal stopping-time models, so there is an opportunity to investigate more complex action spaces. 

\bibliographystyle{unsrt}
\bibliography{refs}

\begin{appendices}
\renewcommand{\thesubsection}{\Alph{subsection}}
\section{Detailed description of the MDP used in the simulation study}\label{appendix: parameters}
We choose to use monthly decision epochs, and, as such, use a standard discount rate $\discount = 0.9975$. The full, unaggregated state space is given by $\states_J = \{0, 1, 2, \ldots, J-1, J, J+1\}$, where states $0$ and $J+1$ are absorbing, and states $1, 2, 3, \ldots, J$ are in increasing order of ``goodness" (i.e., state $x+1$ is preferred to state $x$). The \gls{dm} wishes to avoid the absorbing $0$ ``death" state (e.g., equipment is irreparably damaged or a patient dies before treatment is initiated), and state $J+1$ can be considered an absorbing ``post-intervention" state (i.e., the engine in a car is replaced, or a patient with kidney disease undergoes a kidney transplant). Like \cite{regnier2013state}, we set $J=100$ for all of our simulations. The action space is given by $\actions = \{0, 1\}$, where state $0$ is the ``do nothing" action, and state $1$ is the ``intervention" action. For this study, we use a recursive reward function for $a = 1$ given by $\beta = r(100, 1) = 40$,
$r(\state, 1) = \Big(1-\discount \left(p_{\state 0}^0 - p_{\state+1,0}^0\right)\Big) \cdot r(\state+1, 1)$ for all $\state = 1, \ldots 99$, and $r(101, 1) = 0$.
The value $r(\state, 1)$ is a one-time reward representing the expected remaining lifetime of the system by intervening in state $\state$.
For $a=0$, we set $r(\state, 0) = 1$ for $\state = 1, \ldots 100$, representing the month of system life until the next decision epoch. Additionally, we set $r(0, 1)= r(0, 0) = 0$ since a no-longer-viable system can no longer function, and we set $r(101, 0) = 0$ arbitrarily, as the system cannot feasibly transition to the post-intervention state if the intervention action $1$ is not taken. 

Next, we discuss \cite{regnier2013state}'s \gls{tpm}. $P^0$ has three parameters: $\mu, \gamma, \lambda$ such that $0 \leq \mu, \gamma, \lambda \leq 1$. The parameter $\mu$ represents the probability of remaining in the same state at the next decision epoch, whereas the parameters $\lambda$ and $\gamma$ are such that $\mu\gamma^m$ and $\mu\lambda^m$ denote the probability of declining or improving by $m$  states, respectively at the next decision epoch. Since $0$ is an absorbing state, $p_{00}^0 = 1$ and $p_{0j}^0 = 0 \quad \forall j \neq 0$. \cite{regnier2013state} invoke the condition that
$\lambda =\frac{1-\mu-\gamma)}{(1-\gamma - \mu\gamma)}$ to ensure all rows sum to 1.
To keep the state space consistent between $P^0$ and $P^1$, we add a row and column for state $J+1$ where $p_{J+1,J+1} = 1$ and $p_{i,J+1} = p_{J+1, i} = 0$ for all $i= 0, 1, \ldots, J$. The complete TPM $P^0$ is shown below:

\begin{equation}\label{eq:TPM}
P^0 =\begin{matrix}
0\\
1\\
2\\
\vdots\\
i\\
\vdots\\
J\\
J+1
\end{matrix}
    \begin{pmatrix}
    1&0&0&\dots&0&\dots&0&0&0\\
    \sum_{n=1}^\infty \mu \gamma^n&\mu&\mu\lambda&\dots&\mu\lambda^{i-1}&\dots&\mu\lambda^{J-1}&\sum_{n=J-1}^\infty \mu\lambda^n&0\\
    \sum_{n=2}^\infty \mu \gamma^n&\mu\gamma&\mu&\dots&\mu\lambda^{i-2}&\dots&\mu\lambda^{J-2}&\sum_{n=J-2}^\infty \mu\lambda^n&0\\
    \vdots&\vdots&\vdots&\dots&\vdots&\dots&\vdots&\vdots&\vdots\\
    \sum_{n=i}^\infty \mu \gamma^n&\mu\gamma^{i-1}&\mu\gamma^{i-2}&\dots&\mu&\dots&\mu\lambda^{J-i}&\sum_{n=J-i}^\infty \mu\lambda^n&0\\
    \vdots&\vdots&\vdots&\dots&\vdots&\dots&\vdots&\vdots&\vdots\\
    \sum_{n=J}^\infty \mu \gamma^n&\mu\gamma^{J-1}&\mu\gamma^{J-2}&\dots&\mu\gamma^{J-i}&\dots&\mu\gamma&\sum_{n=0}^\infty \mu\lambda^n&0\\
    0&0&0&\dots&0&\dots&0&0&1
    \end{pmatrix}
\end{equation}

For this study, we use design 1 from \cite{regnier2013state}, which sets $\gamma = 0.85$ and $\mu = 0.1$.

\section{Lumping function}\label{section: appendix lumping}
In this subsection, we introduce the lumping function $s: \states_J \mapsto \states_K$. To determine the number of the original states in $\states_J$ encompassed by each state $h'$ in  $\states_K$ we use the division algorithm. Let $J = qK + r$ such that $q \in \mathbb{Z}$ and $0 \leq r < K$. The quantity $q$ represents the minimum number of original states per aggregated state $h' \in \states_K$. The remainder $r$ represents the number of aggregated states $h' \in \states_K$ containing $q+1$ original states in its definition. 
By the logic that the lower states require greater precision, the lowest $K-r$ states in $\states_K$ contain $q$ original states, while the highest $r$ states in $\states_K$ contain $q+1$ original states. Note that the absorbing states $0$ and $J+1$ remain their own states when aggregating MDP models. 

\section{Alternative TPM Design Result}\label{appendix: alternative TPM}

\begin{figure}[ht]
    \centering
    \begin{subfigure}{.49\textwidth}
    \centering
        \includegraphics*[width=0.95\textwidth]{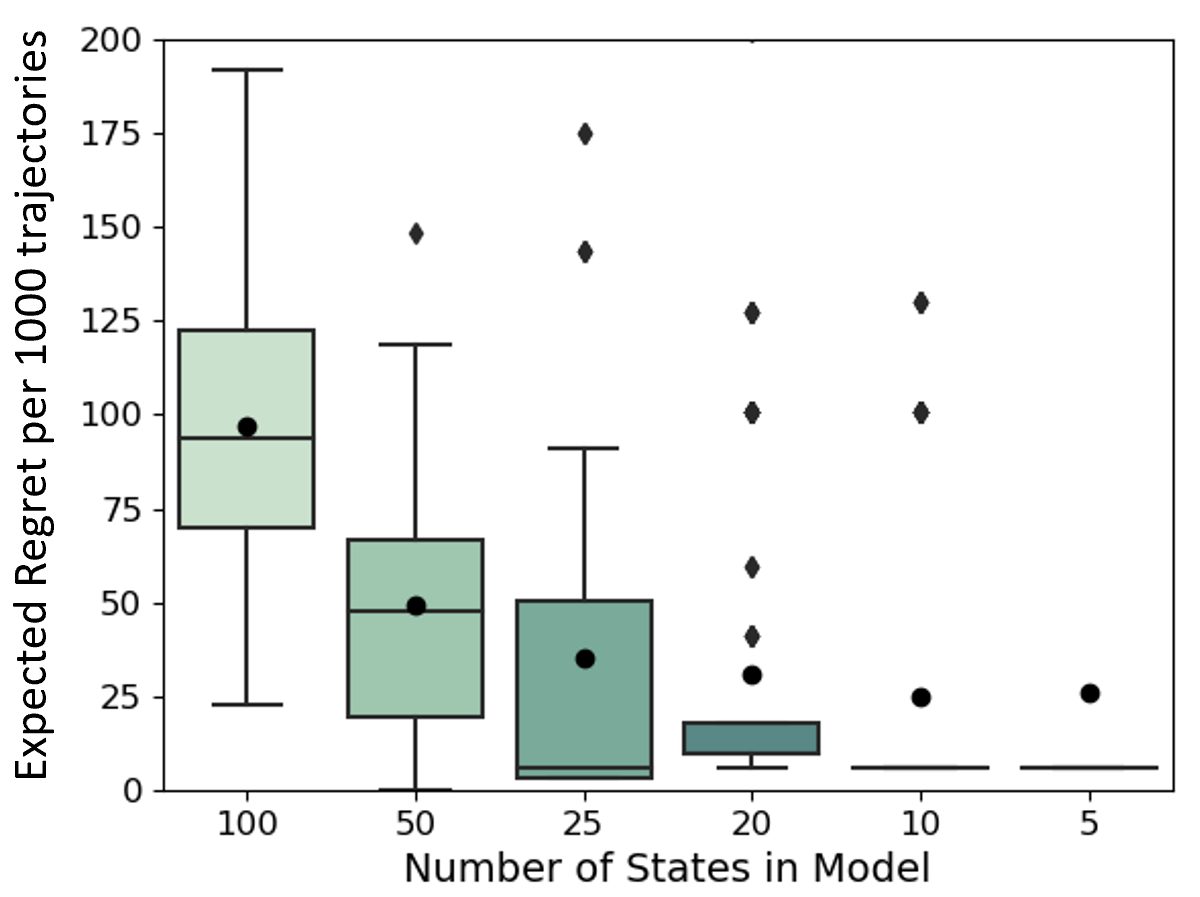}
        \caption{25 system trajectories}
        \label{fig: subfig: regret M=25, Appendix}
    \end{subfigure}
    \begin{subfigure}{.49\textwidth}
    \centering
        \includegraphics*[width=0.95\textwidth]{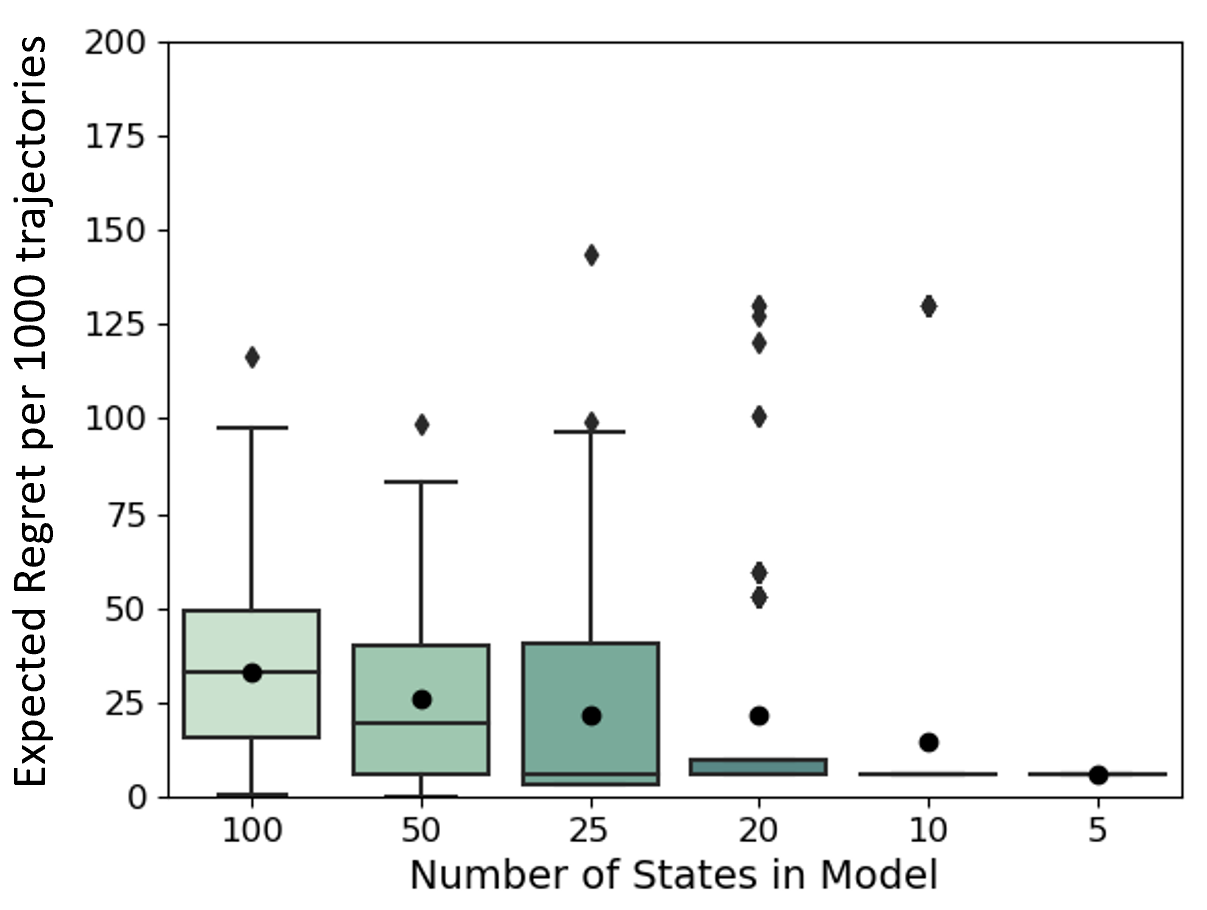}
        \caption{100 system trajectories}
       \label{fig: subfig: regret M=100, Appendix}
    \end{subfigure}
    \caption{The expected regret from making decisions from the optimal policy of an alternative estimated model.}
    \label{fig: expected regret appendix}
\end{figure}

In this Appendix, we give an example of an MDP for which lumping error has minor effects on regret. Let the reward function be the same as the one given in Appendix \ref{appendix: parameters} with the exception that $r(100,T) = 30$. Then, the true threshold is $T_P = 39$. This is an example of a ``lucky" scenario when $T_P$ is extremely close to the maximum value of $L_\state$ for $\state = s(T_P)$.
In Figure \ref{fig: expected regret appendix}, we do not see the U-shaped curve, but rather see that the minimum mean expected regret is given by $K=5$. 


\section{Proofs}\label{appendix: proofs}
\paragraph{Proof of Proposition \ref{prop: Q is DFR}.}
\pf{Suppose not. That is, suppose that $P^a$ has the DFR property, and there is some aggregated matrix $Q^a$ that does not have the DFR property. Because $P^a$ is DFR, it holds that $z_P(i) := \sum_{j = 0}^h p_{ij}^a$
    is non-increasing in $i \in \states_J$ and for all $\state \in \states_J$. If the aggregated matrix $Q^a$ is not DFR, this implies that there exists some $g, k \in \states_K$, such that $
        z_Q(k) = \sum_{l = 0}^g q_{kl}^a < \sum_{l = 0}^g q_{k+1,l}^a = z_Q(k+1)$.
    By the definition of $Q$ in Equation \eqref{eq. p to q}, it follows that for this $g, k$
    \begin{equation}
    \frac{\sum_{l = 0}^g \sum_{i \in L_k} \sum_{j \in L_l} \beta_i p_{ij}^a}{\sum_{i \in L_k} \beta_i} < \frac{\sum_{l = 0}^g \sum_{i \in L_{k+1}} \sum_{j \in L_l} \beta_i p_{ij}^a}{\sum_{i \in L_{k+1}} \beta_i}.
    \label{eq: prop 1 ineq}\end{equation}
    Now, our expression is dependent on the unaggregrated probabilities $p_{ij}^a$ instead of the aggregated probabilities $q_{kl}^a$. Let $m = \max(L_k)$. Because $P$ has the DFR property, it follows that $\sum_{i \in L_k}\sum_{j \in L_l}p_{mj}^a \leq \sum_{i \in L_k}\sum_{j \in L_l}p_{ij}^a$. By definition, we also have that $m + 1 = \min(L_{k+1})$. By a similar logic, because $P^a$ has the DFR property, we have that $\sum_{i \in L_{k+1}}\sum_{j \in L_l}p_{ij}^a\le \sum_{i \in L_{k+1}}\sum_{j \in L_l}p_{m+1, j}^a$. Hence, we can bound \eqref{eq: prop 1 ineq} above and below by the following inequalities:
        \begin{align}
           &&\frac{\sum_{l = 0}^g \sum_{i \in L_k} \sum_{j \in L_l} \beta_i p_{mj}^a}{\sum_{i \in L_k} \beta_i} &< \frac{\sum_{l = 0}^g \sum_{i \in L_{k+1}} \sum_{j \in L_l} \beta_i p_{m+1,j}^a}{\sum_{i \in L_{k+1}} \beta_i}\label{eq: follows DFR}\\
            &\Rightarrow &\frac{\sum_{l = 0}^g \sum_{j \in L_l}p_{mj}^a \sum_{i \in L_k} \beta_i}{\sum_{i \in L_k} \beta_i} &< \frac{\sum_{l = 0}^g \sum_{j \in L_l} p_{m+1,j}^a \sum_{i \in L_{k+1}} \beta_i}{\sum_{i \in L_{k+1}} \beta_i}\label{eq: factor}\\
           &\Rightarrow & \sum_{l = 0}^g \sum_{j \in L_l}p_{mj}^a &< \sum_{l = 0}^g \sum_{j \in L_l}p_{m+1,j}^a\label{eq: reduce}\\
           &\Rightarrow & \sum_{j'=0}^{\max{L_g}}p_{mj'}^a &< \sum_{j'=0}^{\max{L_g}}p_{m+1,j'}^a\label{eq: unlump g},
        \end{align}
    where \eqref{eq: follows DFR} follows from $\transprobs^a$ having the DFR property, \eqref{eq: factor} rearranges the terms in the expression, 
    \eqref{eq: reduce} reduces the like terms in the fraction, and \eqref{eq: unlump g} follows from the disaggregation of the aggregated states $0,..,g\in \states_K$ under the sum. 
    However, this result contradicts the assumption that $P^a$ has the DFR property. Therefore, $Q^a$ must have the DFR property.}\label{proof: Q has DFR property}
\vspace{.5cm}



\hfill
\paragraph{Proof of Proposition \ref{prop: Q gives threshold policy}.}
\pf{It suffices to show that $(\states_K, \actions, Q, \rewards_Q)$ satisfies the three following conditions as adapted from Theorem 3 of \cite{alagoz2004optimal}:
\begin{enumerate}
    \item $Q^0$ has the Decreasing Failure Rate (DFR) property,
    \item $\sum_{l=1}^g q_{kl}^0 \leq \sum_{l=1}^g q_{k - 1, l}^0 \quad \forall k = 2, 3, \ldots, K, \quad g = 1, \ldots, k - 1$, and 
    \item $\frac{r_Q(k, 1) -r_Q(k-1, 1)}{r_Q(k, 1)} \leq \discount(q_{k-1,0}^0 - q_{k 0}^0) \quad \forall k = 2, 3, \ldots K.$ 
\end{enumerate}

Condition 1 is satisfied by Proposition \ref{prop: Q is DFR}. To show condition 2, fix some $k, g$ such that $k \in \{2, \ldots, K\}$ and $g \in \{1, \ldots, k-1\}$. Let $m = \min(L_k)$ and $n = \max(L_{k-1}) = m-1$, and recall that $|L_k| = |L_{k-1}|$ for all $k = 2, 3, 4, \ldots, K$ by the hypothesis. For the left hand side of the inequality in condition 2, we have 

\begin{equation*}
    \sum_{l=1}^g q_{kl}^0 = \sum_{l=1}^g\left(\frac{\sum_{i \in L_k}\sum_{j \in L_l}\beta_i p_{ij}^0}{\sum_{i \in L_k}\beta_{i}}\right) \leq \sum_{l=1}^g\left(\frac{\sum_{i \in L_k}\sum_{j \in L_l}\beta_i p_{mj}^0}{\sum_{i \in L_k}\beta_{i}}\right) = |L_k| \sum_{l = 1}^g \sum_{j \in L_l} p_{mj}^0
\end{equation*}

where the inequality follows from the fact that $P$ has the DFR property. Similarly for the right hand side, we have

\begin{equation*}
\begin{split}
    \sum_{l=1}^g q_{k-1,l}^0 &= \sum_{l=1}^g\left(\frac{\sum_{i \in L_{k-1}}\sum_{j \in L_l}\beta_i p_{ij}^0}{\sum_{i \in L_{k-1}}\beta_{i}}\right) \\
    &\geq \sum_{l=1}^g\left(\frac{\sum_{i \in L_{k-1}}\sum_{j \in L_l}\beta_i p_{nj}^0}{\sum_{i \in L_{k-1}}\beta_{i}}\right)\\
    &= |L_{k-1}| \sum_{l = 1}^g \sum_{j \in L_l} p_{nj}^0.
\end{split}
\end{equation*}

Because $|L_{k-1}| = |L_k|$, and $P$ has the DFR property, 
$$|L_k| \sum_{l = 1}^g \sum_{j \in L_l} p_{mj}^0 \leq |L_{k-1}| \sum_{l = 1}^g \sum_{j \in L_l} p_{nj}^0.$$
Hence, 
\begin{equation*}
    \sum_{l=1}^g q_{kl}^0 \leq |L_k| \sum_{l = 1}^g \sum_{j \in L_l} p_{mj}^0 \leq |L_{k-1}| \sum_{l = 1}^g \sum_{j \in L_l} p_{nj}^0 \leq \sum_{l=1}^g q_{k-1,l}^0,
\end{equation*} as desired. 

Now, we show condition 3 holds. Choose arbitrary $k \in \{2, 3, \ldots, K\}$, and let $m = \min(L_k)$ and $n = \max(L_{k-1})$. Because $(\states_J, \actions, \transprobs, \rewards)$ satisfies the conditions in \cite{alagoz2004optimal}, we know that 
\begin{equation*}
    \frac{r(m, 1)-r(n,1)}{r(m,1)} \leq \alpha(p_{n0}^0-p_{m0}^0).
\end{equation*}
By our assumption, we have 
\begin{equation*}
    \frac{r_Q(k, 1)-r_Q(k-1, 1)}{r_Q(k, 1)} \leq \frac{r(m, 1)-r(n,1)}{r(m,1)}.
\end{equation*} 
Hence, $\frac{r_Q(k, 1)-r_Q(k-1, 1)}{r_Q(k, 1)} \leq \alpha(p_{n0}^0-p_{m0}^0)$. Putting the right hand side in terms of $Q$ gives us:

\begin{align}
\alpha(p_{n0}^0-p_{m0}^0) 
&= \alpha \left(\frac{\sum_{h' \in L_{k-1}}\beta_{h'} p_{n0}^0}{\sum_{h' \in L_{k-1}}\beta_{h'}} - \frac{\sum_{h \in L_{k}}\beta_{h} p_{m0}^0}{\sum_{h \in L_k}\beta_h}\right) \\
&\leq \alpha \left(\frac{\sum_{h' \in L_{k-1}}\beta_{h'} p_{h'0}^0}{\sum_{h' \in L_{k-1}}\beta_{h'}} - \frac{\sum_{h \in L_{k}}\beta_{h} p_{h0}^0}{\sum_{h \in L_k}\beta_h}\right)\label{eq: dfr follows} \\
&= \alpha(q_{k-1,0}^0-q_{k0}^0)
\end{align}

where \eqref{eq: dfr follows} follows from the \gls{dfr} property. Thus, $\frac{r_Q(k, 1)-r_Q(k-1, 1)}{r_Q(k, 1)} \leq \alpha(q_{k-1,0}^0-q_{k0}^0)$, satisfying condition 3. Therefore, by Theorem 3 of \cite{alagoz2004optimal}, the MDP $(\states_K, \actions, Q, \rewards_Q)$ is guaranteed to have a threshold policy that is optimal.}\label{proof: Q gives threshold policy}

\paragraph{Proof of Proposition \ref{prop: Q gives threshold policy in our example}.}
\pf{By Proposition \ref{prop: Q gives threshold policy}, if  $(\states_J, \actions, P, \rewards)$ has a threshold policy that is optimal, and there is some state aggregation function $s: \states_J \mapsto \states_K$ and some reward function $\rewards_{Q_K}$ such that $|L_{k-1}| = |L_{k}|$ and
$\frac{r_Q(k, \, 1) - r_Q(k-1, \, 1)}{r_Q(k \,, 1)} \leq \frac{r(\min(L_{k}), \, 1) -r(\max(L_{k-1}), \, 1)}{r(\min(L_{k}), \, 1)}$, then $(\states_K, \actions, Q, \rewards_Q)$ necessarily has a threshold policy that is optimal. Hence, we need to prove that for all $k \in \states_K \setminus \{0, 1\}$
\begin{enumerate}
    \item $|L_{k-1}| = |L_{k}|$,
    \item $\frac{r_Q(k, \, 1) - r_Q(k-1, \, 1)}{r_Q(k \,, 1)} \leq \frac{r(\min(L_{k}), \, 1) -r(\max(L_{k-1}), \, 1)}{r(\min(L_{k}), \, 1)}$, and
    \item $(\states_J, \actions, P, \rewards)$ used in this study has a threshold policy that is optimal. We can show this using the following three conditions adapted from Theorem 3 of \cite{alagoz2004optimal}:
    \begin{enumerate}
        \item $P^0$ has the Decreasing Failure Rate (DFR) property
        \item $\sum_{l=1}^g p_{\state l}^0 \leq \sum_{l=1}^g p_{\state - 1, l}^0 \quad \forall \state = 1, \ldots, J, \quad g = 1, \ldots, \state - 1$
        \item$\frac{r(\state, 1) r_(\state-1, 1)}{r(\state, 1)} \leq \discount(p_{\state-1,0}^0 - p_{\state 0}^0) \quad \forall \state = 2, 3, \ldots J.$
    \end{enumerate}
\end{enumerate} 
The first criterion is necessarily satisfied in our case study because all values of $K$ evenly divide $J$. The second criterion can be verified algebraically for all values of $k$ using the reward function outlined in Appendix \ref{appendix: parameters}. Proposition \ref{prop: P has DFR property} verifies 3(a) holds.  Since $\gamma \in [0,1]$, $\sum_{k=1}^j p_{\state k}^0 = \sum_{k=1}^j \mu \gamma^{\state - k} \leq \sum_{k=1}^j \mu \gamma^{\state-k-1} = \sum_{k=1}^j p_{\state-1, k}^0$, it follows that 3(b) holds. Finally, 3(c) can be satisfied using basic algebra on the rewards and \gls{tpm} used in this study (Appendix \ref{appendix: parameters}). Thus, by Theorem 3 of \cite{alagoz2004optimal}, the all MDPs denoted by $(\states_K, \actions, Q_K, \rewards_{Q_K})$ used in this study guarantee a threshold policy.}\label{proof: Q gives threshold policy in our example}

\end{appendices}
\end{document}